\documentclass[11pt,a4paper,reqno]{amsart}
\usepackage{calc,amsfonts,amsthm,amscd,epsfig,psfrag,amsmath,amssymb,enumerate,graphicx}
\usepackage[showlabels,sections,floats,textmath,displaymath]{}
\reversemarginpar
\newlength\fullwidth
\numberwithin{equation}{section}

\DeclareMathSymbol{\leqslant}{\mathalpha}{AMSa}{"36} 
\DeclareMathSymbol{\geqslant}{\mathalpha}{AMSa}{"3E} 
\DeclareMathSymbol{\eset}{\mathalpha}{AMSb}{"3F}     
\renewcommand{\leq}{\;\leqslant\;}                   
\renewcommand{\geq}{\;\geqslant\;}                   

\def\1{\ifmmode {1\hskip -3pt \rm{I}} \else {\hbox {$1\hskip -3pt \rm{I}$}}\fi}

\newcommand{\var}{\operatorname{Var}}

\newcommand{\la}{\label } \newcommand{\si}{\sigma } 
\newcommand{\be}{\begin{equation} } \newcommand{\tmix}{T_{\rm mix}} 

\newtheorem{Th}{Theorem}[section]

\newtheorem{Theorem}{Theorem}[section] 
\newtheorem{Le}[Theorem]{Lemma} 
\newtheorem{Proposition}[Theorem]{Proposition} 
 
\newtheorem{remark}[Theorem]{Remark}

\newcommand{\N}{\mathbb N}
\newcommand{\Z}{\mathbb Z}

\newcommand{\cE}{\ensuremath{\mathcal E}} 
 
\newcommand{\cG}{\ensuremath{\mathcal G}}

\newcommand{\cM}{\ensuremath{\mathcal M}}

\newcommand{\bbE}{{\ensuremath{\mathbb E}} }

\newcommand{\bbN}{{\ensuremath{\mathbb N}} } 
 
\newcommand{\bbP}{{\ensuremath{\mathbb P}} } 
 
\newcommand{\bbR}{{\ensuremath{\mathbb R}} }

\newcommand{\bbZ}{{\ensuremath{\mathbb Z}} } 

\let\a=\alpha    \let\d=\delta  \let\e=\varepsilon
 \let\g=\gamma       
      \let\o=\omega      
\let\r=\rho  \let\s=\sigma

\let\O=\Omega      

\def\\{\hfill\break}

\def\thsp{\thinspace}

\def\tthsp{\kern .083333 em}

\def\?{\mskip -10mu}
\newfam\msafam
\newfam\msbfam
\newfam\eufmfam
\def\hexnumber#1{%
  \ifcase#1 0\or 1\or 2\or 3\or 4\or 5\or 6\or 7\or 8\or
  9\or A\or B\or C\or D\or E\or F\fi}
\def\({\left(}
\def\){\right)}

\let\neper=e
\let\ii=i

\def\nep#1{ \neper^{#1}}

\def\tc{\thsp | \thsp}
\let\<=\langle
\let\>=\rangle

\def\gap{\mathop{\rm gap}\nolimits}

\newcommand{\wt}{\widetilde } 

\newcommand{\grad}{\nabla}

\begin{document}

\title[]{Convergence to equilibrium of biased plane partitions}

\author{Pietro Caputo}
\address{Dipartimento di Matematica, Universit\`a Roma Tre,
Largo S.\ Murialdo 1, 00146 Roma, Italia, {\em and} 
UCLA Mathematics Department
Box 951555
Los Angeles, CA 90095-1555.  e--mail: {\tt caputo@mat.uniroma3.it}}
\author {Fabio Martinelli}
\address{Dipartimento di Matematica, Universit\`a Roma Tre,
Largo S.\ Murialdo 1, 00146 Roma, Italia. e--mail: {\tt martin@mat.uniroma3.it}}
\author {Fabio Lucio Toninelli}
\address{CNRS and ENS Lyon, Laboratoire de Physique\\ 46 All\'ee
  d'Italie, 69364 Lyon, France.
e--mail: {\tt fabio-lucio.toninelli@ens-lyon.fr}}

\begin{abstract}
  We study a single-flip dynamics for the monotone surface in $(2+1)$
  dimensions obtained from a boxed plane partition. The surface is
  analyzed as a system of non-intersecting simple paths. When the
  flips have a non-zero bias we prove that there is a positive
  spectral gap uniformly in the boundary conditions and in the size of
  the system. Under the same assumptions, for a system of size $M$,
  the mixing time is shown to be of order $M$ up to logarithmic
  corrections.
\end{abstract}

\maketitle

\thispagestyle{empty}

\section{Introduction, model and results}
Consider a surface in $2+1$ dimensions defined by non-negative integer
heights $\ell_{x,y}$, where $x,y\in\bbZ_+$, such that $\ell_{x,y}\geq
\ell_{x+1,y}$ and $\ell_{x,y}\geq \ell_{x,y+1}$ for all
$x,y\in\bbZ_+$. When $\ell_{x,y}=0$ for all but finitely many $x,y$
this is called a plane partition, the two-dimensional generalization
of an ordinary partition (Young diagram).  When the surface is such
that $\ell_{x,y}\leq c$ and $\ell_{x,y}=0$ when either $x\geq a$ or
$y\geq b$, for some integers $a,b,c$, then it defines a {\sl boxed plane
partition}, or a plane partition in the box $a\times b\times c$.  As we
shall see, a convenient representation of a boxed plane partition is
obtained by considering a system of non-intersecting simple lattice
paths. Other well known equivalent characterizations are the perfect
matchings (or dimers) configurations on a subgraph of the honeycomb
lattice and the lozenge tilings of an 
hexagon.

A continuous time flip dynamics of a plane partition in the box
$a\times b\times c$ is defined as follows: Every point $(x,y)$ in the
rectangle $R_{a,b}=\{0,\dots,a-1\}\times \{0,\dots,b-1\}$ is equipped
with an independent, rate $1$, Poisson clock. When $(x,y)$ rings we
flip an independent $\{0,1\}$ coin $X$; if $X=1$ we replace
$\ell_{x,y}$ by $\ell'_{x,y}=\ell_{x,y}+1$ if allowed; if $X=0$ we
replace $\ell_{x,y}$ by $\ell'_{x,y}=\ell_{x,y}-1$ if allowed. When
the coin is unbiased the flip dynamics converges to the uniform
distribution over plane partitions in the box $a\times b\times c$. The
latter is known to exhibit non trivial limiting shape or ``arctic
circle'' phenomena in the limit of large parameters $a,b,c$, cf.\
\cite{CLP}, \cite{K} and references therein. We refer to
\cite{CK} for the connection with the Wulff crystal in the
low-temperature 3D Ising model.  It is an interesting open problem to
determine the speed of convergence to the uniform equilibrium measure.
We refer to \cite{Wilson} and references therein for the various
polynomial bounds known so far and for the conjectured diffusive
behavior of the spectral gap.

In this paper we shall analyze the biased case. When $X=1$ with
probability $p\neq\frac12$ then the flip dynamics converges to the
probability measure $\mu_\a$ on plane partitions in the box $a\times
b\times c$ such that every configuration $\ell$ has a weight
proportional to $\nep{-2\a {\rm Vol}(\ell)}$, where ${\rm Vol}(\ell)$
stands for the volume $\sum_{x,y} \ell_{x,y}$ under the surface, and
$\nep{-2\a}=p/(1-p)$. By symmetry, $\mu_{\a}$ is equivalent to
$\mu_{-\a}$ and we shall restrict to positive values of $\a$ (or
$p<1/2$). In this model there is no critical value of $\a$ and one has
a localized surface (i.e., ${\rm Vol}(\ell)=O(1)$) for all $\a>0$.
This follows from the fact that the number of plane partitions $\ell$
such that ${\rm Vol}(\ell)=v$ is $\nep{O(v^{2/3})}$, see Section
\ref{prel} for more details.

It has been recently shown that a direct coupling argument allows one to
prove that, if $\a$ is sufficiently large, then uniformly in the size
of the box one has a positive spectral gap and a mixing time 
 of order $M=\max\{a,b,c\}$, see \cite{GPR}. Below we shall prove
that this actually holds for all $\a>0$, up to a logarithmic
correction in the mixing time upper bound. As in \cite{GPR} we use path-coupling arguments
with an exponentially weighted metric. 
However, in the case of small $\a>0$ one of the novelties is that 
these arguments can only be applied
to suitable coarse-grained versions of the process. 
In the proof of the spectral gap estimate we compare the single-flip process to
two auxiliary coarse-grained dynamics. 
Another novel ingredient 
is a roughly deterministic 
description of the process at large scales which allows one to obtain the mixing time estimate.
%

It is worth pointing out that the 
measure $\mu_\a$, $\a>0$, has a natural extension to infinite
boxes, i.e.\ as a measure on plane partitions without any box
constraint. Our spectral gap estimates imply that for all $\a>0$
this extended measure has a positive spectral gap.

\subsection{The model}
We first formulate the model in terms of configurations of 
non-intersecting paths and then describe  
the mapping needed to obtain boxed plane partitions.

Let $k,n,h$ be integers, such that $k,n\geq 1$ and $h\in\{-n,\ldots,n\}$.  We consider
the set $\O^h_{k,n}$ of collections of polymers described as
follows.  Each polymer is a one-dimensional nearest neighbor path of length $n$ which starts at
height $0$ and ends at height $h$, and there are $k$ ordered paths: a
configuration $\eta\in\O^h_{k,n}$ is characterized by integer heights
$\eta^{(j)}_x\in\bbZ$, $j=1,\dots,k$ and $x=0,\dots,n$ satisfying
the constraints:
\begin{gather}\la{mod}
\eta^{(j)}_0=0\,,\quad \eta^{(j)}_{n}=h\,,\quad 
\grad\eta^{(j)}(x):=\eta^{(j)}_{x+1} -
\eta^{(j)}_{x}\in\{-1,+1\}\,,
\;\text{and}\nonumber\\
\eta^{(j)}_x\geq \eta^{(j+1)}_x\,,\quad 1\leq j\leq k-1\,,\; 0\leq x\leq n-1 .
\end{gather} 
The set $\O^h_{k,n}$ is non-empty if $h$ and $n$ have the same parity.
Given $\a>0$, the equilibrium measure $\mu=\mu_{k,n,h}^\a$ is defined
 by \be\la{eqo} \mu(\eta)= \frac{
   \exp\left(\alpha\sum_{i=1}^k\sum_{x=0}^{n}\eta^{(i)}_x
   \right)}{Z}\,,\quad \eta\in\O^h_{k,n}\,,
\end{equation}
where 
$$
Z= Z_{k,n,h}^\a = \sum_{\eta\in \O_{k,n}^h} 
\exp\left(\alpha\sum_{i=1}^k\sum_{x=0}^{n}\eta^{(i)}_x\right)\,
$$
is the normalizing constant.

\begin{figure}[h]
\centerline{
\psfig{file=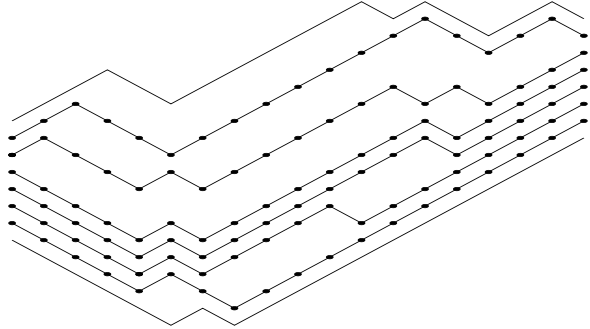,height=1.7in,width=3in}}
\caption{A configuration $\eta\in E_{\xi,\sigma}$. 
Here $k=6,n=18,h=6$.
For graphical convenience $\eta^{(i)}$, 
the $i$-th component of $\eta$, has been shifted by $-i$ units in the vertical direction.
The top path is the ceiling $\xi$ while the bottom path is the floor $\si$ shifted by $-(k+1)$ units.}
\label{k_polymer}
\end{figure}

For every $\xi,\si\in\O^h_{1,n}$, 
we will write simply $\xi\geq\sigma$ 
when $\xi_x\geq \sigma_x$, $ 0\leq x\leq n$.
For a given pair $\xi\geq \si$  we consider 
the subsets
$$
E_{\xi,\si} = \{\eta\in\O^h_{k,n}\,,\;
\xi\geq\eta^{(1)}\, \;\text{ and} \;\,\eta^{(k)}\geq \si\}\,. 
$$
We call $\xi$ and $\si$ the {\em ceiling} and the {\em floor}
respectively, for intuitive reasons, see Figure \ref{k_polymer}. Given a ceiling $\xi$ and a floor
$\si$ ($\xi\geq \si$) we define the equilibrium measure $\mu^{\xi,\si}$ by
conditioning on $E_{\xi,\si}$:
\be\la{eq1}
\mu^{\xi,\si} = \mu\left(\cdot\tc E_{\xi,\si}\right)\,.
\end{equation}
Note that we may write more explicitly
\begin{eqnarray}
\label{eq:moreexp}
  \mu^{\xi,\sigma}(\eta)=\frac{
e^{-2\alpha V^\xi(\eta)}{\bf 1}_{\{\eta\in E_{\xi,\sigma}\}}}
{\hat Z_{\xi,\sigma}}
\end{eqnarray}
where
\begin{eqnarray}
  V^\xi(\eta):=\sum_{1\leq j\leq k}\sum_{1\leq x\leq n-1}\frac{\xi_x-\eta^{(j)}_x}2
\end{eqnarray}
is the ``total volume between $\xi$ and $\eta$'', and
\begin{eqnarray}
\label{eq:hatZ}
  \hat Z_{\xi,\sigma}:=\sum_{\eta\in E_{\xi,\sigma}}
e^{-2\alpha V^\xi(\eta)}.
\end{eqnarray}

We will denote by
$\wedge:=\wedge^{(n,h)}$ the maximal one-polymer configuration in $\O_{1,n}^h$:
\begin{eqnarray}
  \label{eq:wedge}
  \wedge_x:=\left\{
  \begin{array}{lll}
    x&\mbox{for}& 0\leq x\leq (n+h)/2\\
    (n+h)-x&\mbox{for}&(n+h)/2\leq x\leq n.
  \end{array}
\right.
\end{eqnarray}
We will also  use the notation
$\vee:=\vee^{(n,h)} $ for the minimal one-polymer
configuration: 
\begin{eqnarray}
  \label{eq:vee}
  \vee_x:=\left\{
  \begin{array}{lll}
    -x&\mbox{for}& 0\leq x\leq (n-h)/2\\
    (h-n)+x&\mbox{for}&(n-h)/2\leq x\leq n.
  \end{array}
\right.
\end{eqnarray}
Note that $(n+h)/2$ is an integer since $n$ and $h$ have the
same parity, and that $\vee=-\wedge$ if $h=0$.

\subsubsection{From non-intersecting paths to plane partitions}
\label{sec:fromnon}
Suppose first that $\xi = \wedge$ and $\si=\vee$ so that there is no
further constraint on $\eta\in\O^h_{k,n}$. Then there is a bijection
between $\O^h_{k,n}$ and the set of plane partitions in the box
$a\times b\times c$ with $c=k$, $a+b=n$, $b-a=h$ (i.e.\ $a=(n-h)/2$,
$b=(n+h)/2$).  The map is best explained informally as follows. From
the configuration $\eta$ we get the stepped surface $\ell(\eta)$ by
adding layers of height $1$ to the basis rectangle $R_{a,b}$,
according to the paths $\eta^{(j)}$ chosen in reverse order. Namely,
$\eta^{(k)}$ is the first layer. On top of that we put the second
layer $\eta^{(k-1)}$ and so on. 
\begin{figure}[h]
\centerline{
\psfig{file=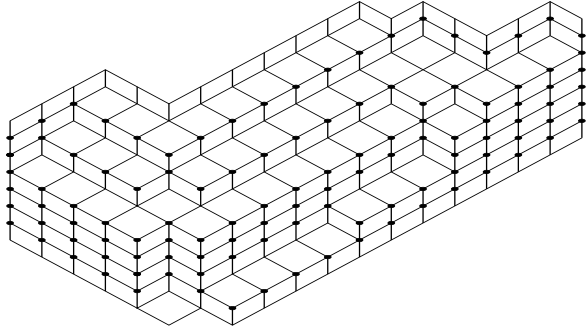,height=1.7in,width=3in}}
\caption{The stepped surface corresponding to the paths in Figure \ref{k_polymer}.}
\label{k_polymer_tiling}
\end{figure}
This defines a bijection, see also Section \ref{sec:polyset} for more details.

For any $\xi\in\O^h_{1,k}$, let $\hat\xi\in\O^{h}_{k,n}$ denote the
configuration such that $\eta^{(i)}=\xi$, $i=1,\dots,k$. Also, let
$\ell(\hat\xi)$ denote the associated stepped surface. Then it is
easily seen that, for any $\xi\geq \si$, the map described above gives
a bijection between the set $E_{\xi,\si}$ and the set of plane
partitions $\ell$ in the box $a\times b\times c$ (as above) such that
$\ell(\hat\xi)\leq \ell\leq \ell(\hat\si)$, see Figure
\ref{k_polymer_tiling}. Note that here inequalities are reversed with
respect to the polymer representation.

\subsection{Results}
\label{sec:results}
The following  heat bath dynamics for configurations of lattice 
paths is easily checked to be equivalent to the flip dynamics discussed in the introduction.
We define the continuous time Markov chain on the set $E_{\xi,\si}$, for a given
pair $\xi,\si\in\O^h_{1,n}$, $\xi\geq \si$, as follows.  At each $(i,x)$,
$i=1,\dots,k$ and $x=1,\dots,n-1$, there is an independent rate $1$
Poisson clock. When $(i,x)$ rings we update the height $\eta^{(i)}_x$
with a new height $\wt \eta^{(i)}_x$ sampled according to the
conditional distribution
$$
\mu_{i,x}^\eta(\cdot):=\mu^{\xi,\si}(\cdot\tc\eta^{(i)}_{x-1}\,,\;\eta^{(i)}_{x+1}
\,,\;\eta^{(i+1)}_x\,,\;\eta^{(i-1)}_{x})
\,,
$$
where $\eta^{(0)}=\xi$ and $\eta^{(k+1)}=\si$.
The Dirichlet form of this process is given by
\be\la{diro}
\cE(f)=\sum_{i=1}^k\sum_{x=1}^{n-1}\mu^{\xi,\si}\left[\var_{i,x}(f)
\right]\,,
\end{equation}
where $\var_{i,x}(f)$ denotes the function 
$$
E_{\xi,\si}\ni\eta\to \var^\eta_{i,x}(f):= \mu_{i,x}^\eta(f^2) -
\mu_{i,x}^\eta(f)^2\,,$$ 
and $f:\O^h_{k,n}\to \bbR$ 
denotes an arbitrary function. The spectral
gap is given by
$$
\gap(\a) = \min_f\,\frac{\cE(f)}{\var(f)}\,,
$$
where $\var(f)=\mu^{\xi,\si}(f^2)-\mu^{\xi,\si}(f)^2$, and the minimum
ranges over all $f:\O^h_{k,n}\to \bbR$ such that $\var(f)\neq 0$. 

\begin{Th}\la{teo_gap}
For any $\a>0$, there exists $c(\a)>0$ such that, uniformly in
$k,n\in\bbN$, $|h|\leq n$  and $\xi\geq\si\in\O^h_{1,n}$:
\be\la{unifo0}
\gap(\a)\geq c(\a)\,.
\end{equation} 
\end{Th}
The estimate in Theorem \ref{teo_gap}
is already known to hold when $k=1$ (see \cite[Th.\ 4.3]{CM} or, for an alternative proof, \cite{GPR}). 

\medskip

Our second result concerns the {\em mixing time} of the Markov chain which,
we recall, is defined as
\begin{eqnarray}
  \label{eq:Tmix}
  \tmix =\inf\{t>0:\max_{\eta\in\Omega^h_{k,n}}\|P_t(\eta,\cdot)-\mu\|_{\rm var}
\leq 1/(2e)\},
\end{eqnarray}
with $\|\cdot\|_{\rm var}$ denoting the total variation norm:
\begin{eqnarray}
  \label{eq:totvar}
  \|\nu-\nu'\|_{\rm var}:=\frac12\sum_{\eta\in \Omega^h_{k,n}}|\nu(\eta)-
\nu'(\eta)|.
\end{eqnarray}
$P_t(\eta,\cdot)$ is the law, at time $t$, of the Markov chain 
started from $\eta$ at time zero.

\smallskip

Theorem \ref{teo_gap} implies that the mixing time
of the Markov chain defined above in the case $2k=n$, $h=0$, $\xi=\wedge$,
$\sigma=\vee$ is
$O(n^3)$. This is a simple consequence of the well-known inequality
\begin{eqnarray}
\label{eq:tmixgap}
 \tmix \leq \gap^{-1}\left(1-\log\left(\min_\eta\mu(\eta)\right)\right),
\end{eqnarray}
(see also Lemma \ref{th:fnraf} below). We can however prove:
\begin{Th}
  \label{th:tmix1}
Let $M=\max(n,k)$. 
  For every $\alpha>0$ there exists
  $C(\alpha)<\infty$ such that uniformly in the choice of ceiling $\xi$ and floor $ \sigma$, 
$\tmix\leq C(\alpha)\, M (\log
  M)^{6}$.
\end{Th}
We mention that the conjectured behavior is $O(M)$, without
logarithmic corrections. 
This bound is known to hold if $\a$ is sufficiently large \cite{GPR}.

For simplicity we have stated these results for a positive constant bias $\a>0$, but there is no difficulty to extend them to the case of a non-homogeneous bias 
$\a_{x,y}$ on each column of the stepped 
surface $\ell_{x,y}$, 
provided that there exists $\a_0>0$ such that $\a_{x,y}\geq \a_0$ for all $x,y$.
Alternatively, one could place a non-homogeneous bias $\a_{(i,x)}$ associated to each polymer $i$ and position $x$ in the definition of the measure (\ref{eqo}).

The rest of the paper is organized as follows.
In Section \ref{prel} we provide some preliminaries and prove a couple of key equilibrium 
estimates to be used in the proof of the main theorems. The latter is given in 
Section \ref{gap} (Theorem \ref{teo_gap}) and Section \ref{mixing} (Theorem \ref{th:tmix1}).

\section{Preliminaries}\la{prel}

\subsection{Particles and vacancies}
Each polymer $\eta^{(j)}$ can be characterized by
  the positions of its positive increments, also called {\em
  particles}. More precisely, for every $1\leq j\leq k$, 
let $x_i^{(j)}=x_i^{(j)}(\eta)$, $i\geq 1$ denote the position of the $i$-th positive
increment in the $j$-th polymer, defined recursively by:
\begin{gather*}
x_1^{(j)}= \min\{x\in \{0,\ldots,n-1\}\,:\; \grad\eta^{(j)}(x)=+1\}\,,\;\dots\,
\\
 x_{\ell+1}^{(j)}= \min\{x> x_\ell^{(j)}\,:\; \grad\eta^{(j)}(x)=+1\}\,.
\end{gather*} 
Note that, given $h$, the $k$ polymers in the configuration
$\eta\in\O^h_{k,n}$ all have the same number $N= N(h,n)=(n+h)/2$ of particles
(they all have the same length $n$, the same starting point $0$ and
the same end-point $h$). Observe that, because of the order
constraint, particles obey the following relations: $x_i^{(j)} \leq
x_i^{(j+1)}$, $i=1,\dots,N$ and $j=1,\dots,k-1$.
We often write ${\bf x}$ or ${\bf x}(\eta)$ for the collection of particle positions of a given 
configuration $\eta$.

The set of {\sl vacancies} for the polymer $\eta^{(j)}$ 
is defined as the set of points in $\{0,\ldots,n-1\}$ which do not
contain particles.  Of course, the number of particles plus the number
of vacancies for $\eta\in\O^h_{1,n}$ equals $n$.

\subsection{Monotonicity}\label{sec:monoton}

Trajectories of the Markov chain corresponding to distinct initial conditions and/or distinct boundary constraints can be realized on the same probability space by a standard coupling argument.
This is a straightforward generalization of the argument for a single polymer, see \cite[Section 2]{CMT}. It follows that the Markov chain enjoys the following useful monotonicity property. If $\eta^{\xi,\si}(t;\zeta)$ denotes the evolution of the surface with ceiling $\xi$ and floor $\si$ at time $t$ and with starting configuration $\zeta$ at time $0$, then almost surely one has
\be\la{mono}
\eta^{\xi,\si}(t;\zeta)\geq \eta^{\xi',\si'}(t;\zeta')\,,
\end{equation}
whenever $\xi\geq\xi'$, $\si\geq\si'$ and $\zeta\geq \zeta'$.
Here for two systems of polymers $\eta,\zeta\in\O^{h}_{k,n}$ we use the convention that $\eta\geq \zeta$ means $\eta^{(i)}\geq \zeta^{(i)}$ for all $i$. 

Let $\bbE$ denote expectation with respect to this global coupling
$\bbP$. Using the notation $\bbE[f(\eta^{\xi,\si}(t;\zeta))] = 
P_t f(\zeta)$, $f:\O^{h}_{k,n}\to\bbR$, then (\ref{mono}) implies
that for every fixed $t\geq 0$, the function $P_t f$ is
increasing whenever $f$ is increasing, where a function $f$ is called
increasing if $f(\eta)\geq f(\zeta)$ for any $\eta,\zeta$ 
such that $\eta\geq \zeta$. 
Useful inequalities for the equilibrium measures can be
derived from this. For instance, taking 
the limit $t\to\infty$ in (\ref{mono}) 
yields the inequality 
\be\la{monoeq}
\mu^{\xi,\si}(f)\geq
\mu^{\xi',\si'}(f)\,, 
\end{equation}
for any increasing $f$ and any 
$\xi\geq \xi', \si\geq\si'$. 
We will often use one form or another of the inequality (\ref{monoeq}) without explicit reference.

\subsection{Tightness of the excess volume and decay of correlation}
Here we prove some equilibrium results concerning the
exponential decay of spatial correlations, and the exponential
tightness of the ``excess volume'' $V^\xi(\eta)$. While our main aim
is to provide the necessary tools for the proof of Theorems
\ref{teo_gap} and \ref{th:tmix1}, such results may be of independent
interest.

We start with a basic estimate for the case where $h=0$, $n\in2\N$ and the
ceiling is the maximal possible configuration, $\wedge$. 

\begin{Le}
  \label{lemma_basico}Let $h=0$.  For every $\alpha>0$ there exists
  $c_1(\alpha)>0$ such that, uniformly in $n\in2\N$, $k\in\N$ and
  in the configuration $\sigma\in\O^0_{1,n}$ of the floor,
\begin{eqnarray}
\label{eq:cl1}
   \mu^{\wedge,\sigma}(V^\wedge(\eta)\geq i)\leq e^{-c_1(\alpha)i}
\end{eqnarray}
for every $i\in\N$. In particular, there exists $p(\alpha)>0$ such that
\begin{eqnarray}
\label{eq:cl2}
  \mu^{\wedge,\sigma}\left(\eta^{(j)}=\wedge,\;\forall\, 1\leq j\leq k\right)
\geq p(\alpha).
\end{eqnarray}
\end{Le}

\proof By monotonicity, we can upper bound the probability in
\eqref{eq:cl1} replacing $\sigma$ with $\vee$: 
\begin{eqnarray}
\label{eq:Zwv}
   \mu^{\wedge,\sigma}(V^\wedge(\eta)\geq i)\leq
    \mu^{\wedge,\vee}(V^\wedge(\eta)\geq i).
\end{eqnarray}
Always by monotonicity, the right-hand side of \eqref{eq:Zwv} is 
non-decreasing in $n$ and $k$. Therefore, 
\begin{eqnarray}
\label{eq:Nv}
  \mu^{\wedge,\sigma}(V^\wedge(\eta)\geq i)\leq
  \frac{\sum_{v\geq i}e^{-2\alpha v}\mathcal N(v)}
  {\sum_{v\geq 0}e^{-2\alpha v}\mathcal N(v)}
\end{eqnarray}
where $\mathcal N(v)$ is the number of plane partitions of volume $v$,
and the right-hand side is just the limit $n\to\infty,k\to\infty$ of
the right-hand side of \eqref{eq:Zwv}.  
The
dependence on $\sigma,k$ and $n$ has then disappeared.  Since
\begin{eqnarray}
\label{eq:ubNv}
\mathcal N(v) \stackrel{v\to\infty}\sim 
\frac{a_1}{v^{25/36}}\exp\left(a_2 v^{2/3}\right)
\end{eqnarray}
for some (explicit) positive constants $a_1,a_2$ 
 \cite{cf:Wright}, one obtains
immediately \eqref{eq:cl1}.  Eq. \eqref{eq:cl2} is obtained from \eqref{eq:cl1}
just taking $i=1$ and $p(\alpha)=1-\exp(-c_1(\alpha))$. \qed

\medskip


Even if the ceiling $\xi$ does not coincide with $\wedge$, 
it is intuitive that the polymer $\eta^{(1)}$ 
 gets more and more squeezed to $\xi$ when $k$ becomes
large. This is one of the implications of the next result:
\begin{Le}
\label{lemma_spiaccicato} For every $\alpha>0$ there exists $c_2(\alpha)>0$
such that the following holds.
  Uniformly in $\xi\geq\sigma\in \O^h_{1,n}$, in $n,k\in \N$, 
$|h|\leq n$ and $0\leq x\leq n$,
\begin{eqnarray}
\label{eq:spiacc1}
  \mu^{\xi,\sigma}(\eta^{(j)}_x\neq\xi_x)\leq e^{-c_2(\alpha)(k-j+1)}.
\end{eqnarray}
Moreover, for every $0<a<b<n$, 
\begin{eqnarray}
\label{eq:spiacc2}
  \mu^{\xi,\sigma}\left(\eta^{(k)}_x\ne \xi_x\,,\;\forall a\leq x\leq b\right)\leq
e^{-c_2(\alpha)(b-a)}.
\end{eqnarray}
\end{Le}

\begin{figure}[!h]
\centerline{
\psfrag{x}{$x$}
\psfrag{y}{$y$}
\psfrag{x}{$x$}
\psfrag{0}{$0$}
\psfig{file=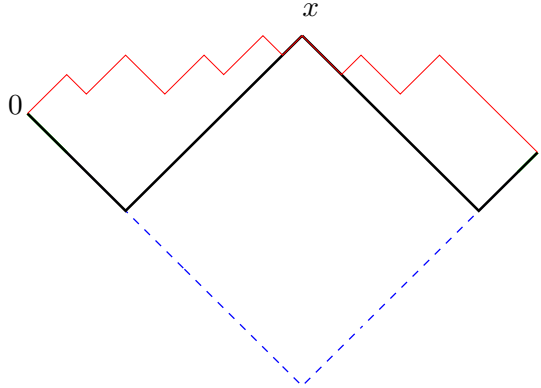,height=2in}}
\caption{Graphical construction of $\tilde\wedge$. The thin full line
  denotes $\xi$, the dashed line $\vee$ and the thick full line is
  $\tilde\wedge$, for a given value of $x$.  The floor $\sigma$ is not
  drawn, since it has no influence on the construction of
  $\tilde\wedge$.  }
\label{fig:Vtilde}
\end{figure}

\proof 
For a fixed $1\leq x<n$, let
$\tilde \wedge\in\Omega_{1,n}^h  $ be 
defined as follows (see Fig. \ref{fig:Vtilde}): 
\begin{eqnarray}
  \tilde\wedge_y:=\max(\xi_x-|y-x|,\vee_y)\,,\;0\leq y\leq n.
\end{eqnarray}
Note that $\xi\geq \tilde\wedge$ and of course $\sigma\geq \vee$.
Then, by  monotonicity and the fact that $\xi_x=\tilde \wedge_x$ we have
\begin{eqnarray}
  \mu^{\xi,\sigma}(\eta^{(j)}_x<\xi_x)\leq
  \mu^{\tilde\wedge,\vee}(\eta^{(j)}_x<\tilde\wedge_x).
\end{eqnarray}
Now, $\eta^{(j)}_x< \tilde\wedge_x$ implies $V^{\tilde\wedge}(\eta)\geq
k-j+1$. As in Eq. \eqref{eq:Nv}, from monotonicity it follows that
\begin{eqnarray}
\label{eq:Nv2}
  \mu^{\xi,\sigma}(\eta^{(j)}_x<\xi_x)\leq
  \frac{\sum_{v\geq (k-j+1)}e^{-2\alpha v}\mathcal N(v)}
  {\sum_{v\geq 0}e^{-2\alpha v}\mathcal N(v)}
\end{eqnarray}
 and \eqref{eq:spiacc1} follows from \eqref{eq:ubNv}.

 Next, we prove \eqref{eq:spiacc2}. Denote by $x^{(0)}_r$, $1\leq r\leq
 N(n,h)$ the positions of the particles of the ceiling $\xi$, and let
 $\mathcal I_{a,b}:=\{x^{(0)}_{i},x^{(0)}_{i+1},\ldots,
 x^{(0)}_{i+m}\}$ the set of those particle positions which are contained in
 the interval $\{a,\ldots,b-1\}$ (the cardinality $m+1$ of $\mathcal
 I_{a,b}$ does not exceed $(b-a)$ and can be zero). We use also the
 notation $x^{(k+1)}_r$, $1\leq r\leq N(n,h)$ to denote the positions of
 the particles of the floor $\sigma$. The event in the left-hand side
 of \eqref{eq:spiacc2} implies that $x^{(k)}_{r}>x^{(0)}_r$ for every
 $i\leq r\leq i+m$. On the other hand, from \eqref{eq:spiacc1} we know
 that the event $\{x^{(j)}_{i+m}=x^{(0)}_{i+m} \,\forall\; 1\leq j\leq
 k\}=\{x^{(k)}_{i+m}=x^{(0)}_{i+m}\}$ has probability at least
  $1-\exp(-c_2(\alpha))$, uniformly in
 all the parameters.  Assume that this event is not realized. In this case,
 the probability that $\{x^{(j)}_{i+{m-1}}=x^{(0)}_{i+{m-1}} \,\forall\;
 1\leq j\leq k\}=\{x^{(k)}_{i+m-1}=x^{(0)}_{i+m-1}\}$ is again lower bounded by $1-\exp(-c_2(\alpha))$:
 indeed, by monotonicity it is sufficient to consider the case where
 $x^{(j)}_{i+m}=x^{(k+1)}_{i+m}$ for every $1\leq j\leq k$, and to apply
 once more \eqref{eq:spiacc1}. Iterating this procedure, we see that
 the left-hand side of \eqref{eq:spiacc2} is lower-bounded by
\begin{eqnarray}\la{iab}
  e^{-c_2(\alpha)|\mathcal I_{a,b}|}.
\end{eqnarray}
One can then repeat the argument with the vacancies replacing the 
particles. The argument is the same except 
that vacancies have to be matched from left to right (while particles have been matched from right to left). Since the number of vacancies plus the number of particles
in $\{a,\ldots,b\}$ equals $(b-a)$, one obtains immediately 
\eqref{eq:spiacc2} (modulo redefining $c_2(\alpha)$).

\qed

\section{Proof of Theorem \ref{teo_gap}}\la{gap}
To prove the spectral gap estimate we shall use a three-fold decomposition 
that can be roughly described as follows. The first step, carried out in Section \ref{fulll},
allows to reduce the original process to a process with a bounded (i.e.\ independent of $n,k$) 
number $s$ of polymers. This is achieved by estimating the spectral gap of 
a new process whose elementary moves consist in the updatings 
of sets of $s$ full-length adjacent polymers. 
The second step, see Section \ref{particles} below, allows to 
further reduce the problem to a process 
where each of the $s$ polymers has a bounded
(i.e.\ independent of $n,k$) number $r$ of particles. This is obtained by analyzing a 
coarse-grained dynamics where we update blocks of particles. 
The last step is a rough estimate for a system with $s$ polymers and 
$r$ particles, see Section \ref{bounded}.  Finally, in 
Section \ref{together} we prove the theorem by gathering 
all the pieces together.

\subsection{Particle block--dynamics}\la{particles}
%
%
%
%
%


Let $\nu^{{\bf x}}_{j,m}$, for $j<m$, denote the equilibrium measure
$\mu^{\xi,\si}$ on $E_{\xi,\si}$
conditioned on the $\s$-algebra generated by the particle positions
\be\la{partcond2}
\{x_{v}^{(u)}\,,\;u=1,\dots,k\}\,,\quad v \leq  j \,,\;
\text{and}\;
v\geq m\,.
\end{equation}
As a convention, if $v<1$ then we set (deterministically)
$x_v^{(j)}=0$ for all $j$. Similarly, if $v> N$, we set $x_v^{(j)}=N$.
Recall that $N=N(h,n)$ denotes the total number of particles.

\begin{figure}[h]
\centerline{
\psfig{file=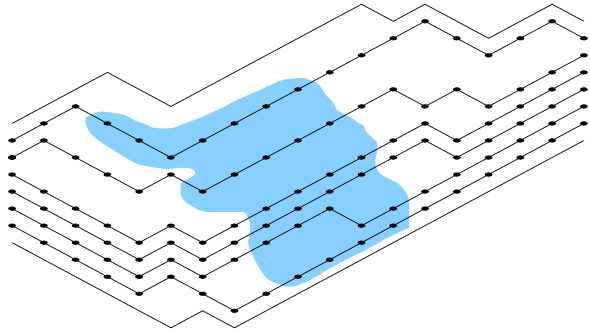,height=1.7in,width=3in}}
\caption{An illustration of the definition of the measure $\r^{\bf x}_{i,\ell}$ for the paths 
from Figure \ref{k_polymer}. 
Here $\ell=2$, $i=5$. The shaded region is the portion of paths distributed according to 
$\r^{\bf x}_{i,\ell}$, while the rest of the configuration is frozen.}
\label{rhos}
\end{figure}

It will be convenient to have the following alternative notation.  For
a fixed integer $\ell$ we define the measures $\r^{{\bf
    x}}_{i,\ell}:=\nu^{{\bf x}}_{i-\ell-1,i+\ell+1}$.  Once the values
of (\ref{partcond2}) are given, $\r^{{\bf x}}_{i,\ell}$ is a
distribution on the portion of paths $\eta^{(j)}$ 
in the segments $\{a^{(j)}+1,\dots,b^{(j)}-1\}$, $j=1,\dots,k$, where
$a^{(j)}:=x^{(j)}_{i-\ell-1}$, $b^{(j)}:=x^{(j)}_{i+\ell+1}$, see
Figure \ref{rhos}. 

We will study the following Markov chain with state space
$E_{\xi,\si}$,
for a given pair $\xi,\si\in\O^h_{1,n}$, with $\xi\geq\si$.
This auxiliary process, which we call the particle block-dynamics, is reversible
w.r.t.\ the equilibrium measure $\mu^{\xi,\si}$.

We have $N(h,n)$ independent Poisson clocks with parameter
$1$. When the $i$--th clock rings we consider 
the current configuration $\eta$, and update the portion of paths
$\eta^{(j)}$ in the segments $\{a^{(j)}+1,\dots,b^{(j)}-1\}$,
$j=1,\dots,k$
with a sample from the conditional distribution $\r^{\bf x}_{i,\ell}$.
The rest of the configuration $\eta$
is left unchanged. In other words, we are removing from the system 
all particles 
at positions $x_{i-\ell}^{(j)},\dots,x_{i+\ell}^{(j)}$, $j=1,\dots,k$
and
we are replacing them by a sample from $\r^{\bf x}_{i,\ell}$.

\smallskip

The generator of this process can be written as 
\be\la{geng}
\cG f = \sum_{i=1}^N \left[\r^{{\bf
    x}}_{i,\ell}(f) - f\right]\,,
\end{equation}
where $f$ denotes a function $f:E_{\xi,\si}\to\bbR$ and $\r^{{\bf
    x}}_{i,\ell}(f)$
is the function $\eta\to\int f(\zeta)\r^{\bf x}_{i,\ell}(d\zeta)$ for ${\bf x}={\bf x}(\eta)$.
Since $\r^{{\bf
    x}}_{i,\ell}$ are conditional expectations we see that 
$\mu^{\xi,\si}\left(f(\r^{{\bf
    x}}_{i,\ell}(f) - f)\right) = 
\mu^{\xi,\si}\left((\r^{{\bf
    x}}_{i,\ell}(f))^2 - \r^{{\bf
    x}}_{i,\ell}(f^2)\right)$, so
that the Dirichlet form of the process is
\be\la{dirrho}
-\mu^{\xi,\si}\left(f\cG f\right) = \sum_{i=1}^N 
\mu^{\xi,\si}\left(\var_{\r^{{\bf
    x}}_{i,\ell}}(f)\right)\,.
\end{equation}
In particular, $\cG$ is self-adjoint 
in $L^2(\mu^{\xi,\si})$. 
Let $\gap(\cG)$ denote the spectral gap of this process, which 
of course depends on the choice of $\ell\in\N$.

\begin{Proposition}\la{pros}
For any $\a>0$ and any $k\in\bbN$, there exists $\ell=\ell(\a,k)$
such that uniformly in $h$, $n$, and $\xi,\si\in\O^h_{1,n}$ we have
$$
\gap(\cG)\geq 1\,.
$$
\end{Proposition}

To prove Proposition \ref{pros} we use a coupling argument.  Consider
two evolutions $\eta(t),\eta'(t)$, $t\geq 0$, of the Markov chain
described above, with initial conditions $\eta$ and $\eta'$
respectively, where $\eta,\eta'\in E_{\xi,\si}$. 
A well-known argument (see e.g.\ Proposition 3 in
\cite{Wilson}) shows that, for any coupling $\bbP$ of the two
evolutions 
\be \gap(\cG) \geq - \liminf_{t\to\infty} \frac1{t}\log
\left[\max_{\eta,\eta'}\,\bbP\left(\eta(t)\neq
    \eta'(t)\right)\right]\,.  \la{gapcoup}
\end{equation}
Therefore, to prove Proposition
\ref{pros} it is sufficient to establish 
Lemma \ref{proel} below.

\begin{Le}\la{proel}
For every $\a>0$ and $k\in\bbN$ there exist $\ell\in\N,\g>0$ and 
a coupling $\bbP$ of $(\eta(t),\eta'(t))$
such that, uniformly in the starting configurations and uniformly in
the parameters $h,n$ and $\xi,\si\in\O_{1,n}^h$ 
$$
\bbP(\eta(t) \neq \eta'(t))\leq k N(h,n)\nep{\g\,n}\,\nep{-t}\,,\quad t\geq 0\,.
$$
\end{Le}
\proof
Consider the distance
\be\la{dis1}
d_\g(\eta,\eta') = \sum_{i=1}^N\sum_{j=1}^k
\nep{-\g\,i} 1_{\{x_i^{(j)}\neq
  y_i^{(j)}\}}\,,
\quad \eta,\eta'\in E_{\xi,\si}\,,
\end{equation}
where $x_i^{(j)}, y_i^{(j)}$ denote the positions of the $i$-th
particle of the $j$-th polymer in the configurations $\eta$ and
$\eta'$ respectively, and $\g>0$ is to be determined later. 
Note that the minimal non-zero value of $d_\g(\cdot,\cdot)$
is $\nep{-\g n}$. Therefore, by Markov's inequality
$$
\bbP(\eta(t) \neq \eta'(t))\leq \nep{\g\,n}\,\bbE(d_\g(\eta(t),\eta'(t)))\,.
$$
We need to show that we can define a coupling such that
$$\bbE(d_\g(\eta(t),\eta'(t)))\leq k N(h,n)\nep{-t}\,.$$ 
From path coupling (\cite{BD}), see Lemma \ref{path} below for the details, it will be sufficient 
to exhibit a Markovian coupling such that 
\be\la{couple1} [\wt\cG \,d_\g\,](\eta,\eta'):=
\left.\frac{d}{dt}\,\bbE(d_\g(\eta(t),\eta'(t)))\,\right|_{t=0^+}\leq
-\,d_\g(\eta, \eta')
\,,
\end{equation}
for all pairs $\eta,\eta'$ satisfying $d_0(\eta,\eta')=1$,
i.e.\ when there is only one discrepancy in the particle configurations.
In this case, $d_\g(\eta,\eta')=\nep{-\g\,i}$ where $i$ is the label
of the discrepancy. 

Next, we specify the Markovian coupling, and prove that it satisfies (\ref{couple1}). 
If we use the same Poisson clocks for the two evolutions, the infinitesimal 
generator $\wt \cG$ of the coupled dynamics can be written as
\begin{eqnarray}
\label{upto}
\wt\cG = \sum_{i=1}^N (E_{i,\ell} - 1)\,,  
\end{eqnarray}
where $E_{i,\ell}$ denotes a (not yet specified) 
coupling of the local equilibria $\r^{{\bf
    x}}_{i,\ell}$ 
of the block of $2\ell+1$ particles
around the $i$-th particle for the $k$ polymers, cf.\ (\ref{geng}). 

To prove (\ref{couple1}), we may choose the
coupling in such a way that if $\eta,\eta'$ 
have a single discrepancy at a given particle label $i$ and at a given
polymer label $j$, then 
$[E_{v,\ell}d_\g](\eta,\eta') = 0$ for all $v$ such
that  $v-\ell\leq i\leq v+\ell$. 
By construction, there are at least $\ell$ such blocks for any fixed $i$
(this is the case if e.g.\ $i=1$ or $i=N$).

Let $E^-_i=E_{i-\ell-1,\ell}$ (respectively, $E^+_i=
E_{i+\ell+1,\ell}$) denote the coupling corresponding to the block of
$2\ell+1$ particles just to the left (resp.\ just to the right) of
particle $i$.  Note that if e.g.\ $i\leq \ell+1$ then there is no block
just to the left of $i$ and we may set $E^-_i=1$ for such $i$. Similarly, 
if $i\geq N-\ell$ then there is no block just to the right of $i$ and
we can set $E^+_i=1$ in this case.  Since all other blocks give a
trivial contribution to (\ref{couple1}) we see that 
\begin{align}\la{see}
&[\wt\cG \,d_\g\,](\eta,\eta') \\
&\quad
\leq -\ell\,d_\g(\eta,\eta')
+ [(E^-_i-1)\,d_\g\,](\eta,\eta') +
[(E^+_i-1)\,d_\g\,](\eta,\eta')\nonumber\,.
\end{align}
Recall that $d_\g(\eta,\eta') = \nep{-\g\,i}$. 
Then we can estimate
\begin{align}
[E^+_i\,d_\g\,](\eta,\eta') &=
\sum_{u=1}^k\sum_{v=i+1}^{i+2\ell+1} 
\nep{-\g\,v}\,
E_i^+\left[1_{\{x_v^{(u)}\neq y_v^{(u)}\}}\right] \nonumber\\
&\leq k\,\nep{-\g\,i}
\sum_{v=1}^{\infty} \nep{-\g\,v} = 
\frac{k\,\nep{-\g}}{1-\nep{-\g}}\,d_\g(\eta,\eta')\,,
\la{seeo}
\end{align}
where 
we have bounded by $1$ the probability of
a discrepancy.

On the other hand, denoting by $F_v$ the event that there exists 
$u\in\{1,\dots,k\}$ such that $ x_{i-v}^{(u)}\neq y_{i-v}^{(u)}$, we have
\begin{align*}
[E^-_i\,d_\g\,](\eta,\eta') 
&= \sum_{u=1}^k
\sum_{v=i-2\ell-1}^{i-1} \nep{-\g\,v}\,
E_i^- \left[1_{\{x_v^{(u)}\neq y_v^{(u)}\}}\right]\\
&\leq 
k\,\nep{-\g\,i} 
\sum_{v=1}^{2\ell+1} \nep{\g\,v}
E_i^-\left[1_{F_v}\right]  = k\,d_\g (\eta,\eta')
\sum_{v=1}^{2\ell+1} \nep{\g\,v}
E_i^-\left[1_{F_v}\right] \,,
\end{align*}
with the convention that $1_{F_v}=0$ if
$i-v<1$.

Below, see discussion after (\ref{lemmo1}),  we prove 
that for a suitable choice of the coupling 
there exists $c=c(\a,k)>0$ independent of $\gamma$ such that that 
\be\la{claa}
E_i^-\left[1_{F_v}\right]
\leq \nep{-c\,v}\,,\quad\,1\leq v\leq 2\ell+1\,.
\end{equation} 
If we assume this estimate,  
from (\ref{see}), (\ref{seeo}) we conclude
that if e.g.\ $\g=c/2$ then (\ref{couple1}) follows for $\ell$
sufficiently large (depending on $k$ and $\a$). In particular, the proof of Lemma \ref{proel}
will be completed once we prove (\ref{claa}).

\medskip

We turn to the proof of (\ref{claa}). 
Let ${\bf x},{\bf y}$ denote the collections
$\{x_{v}^{(u)}\},\{y_{v}^{(u)}\}$ 
of all positions
of particles of two configurations $\eta,\eta'\in E_{\xi,\si}$
and consider the associated probability measures $\nu^{{\bf x}}_{j,m},\nu^{{\bf y}}_{j,m}$ 
defined in (\ref{partcond2}).   
For a fixed pair of integers $j,m$, let 
$\wt \nu = \wt \nu^{j,m}$ be the independent
coupling of 
$\nu^{{\bf x}}_{j-m-1,j}$ and $\nu^{{\bf y}}_{j-m-1,j}$. That is, 
we are freezing all particles labeled $i\leq j-m-1$ or $i\geq j$ 
and we sample the $m$ particles labeled $i=j-m,\dots,j-1$ 
according to the independent coupling of $\nu^{{\bf x}}_{j-m-1,j}$ 
and $\nu^{{\bf y}}_{j-m-1,j}$.
We say that ${\bf x},{\bf y}$
agree up to $j-m-1$ if $x_{v}^{(u)}=y_{v}^{(u)}$, for all
$u=1,\dots,k$ and for all $v\leq j-m-1$. 
%

\medskip

We claim that 
there exists $\e>0$ depending only on $\a$ and $k$
such that, 
if ${\bf x},{\bf y}$
agree up to $j-m-1$, then
\be\la{lemmo1}
\wt \nu\left(
x_{j-1}^{(u)}=y_{j-1}^{(u)}\,,\;\text{ for all}\;
\;u=1,\dots,k 
\right)\geq \e\,.
\end{equation}

\medskip

Let us first show that (\ref{lemmo1}) implies (\ref{claa}).
First of all, let us sample $x_{i-1}^{(u)},y_{i-1}^{(u)}$,
$u=1,\dots,k$, using $\wt \nu$ with $j=i$
and (supposing for simplicity $i>2\ell+1$) $m=2\ell+1$
(recall that $i$ is the index appearing in the proof of Lemma \ref{proel}).
Then (\ref{lemmo1}) implies that there is a full matching 
$x_{i-1}^{(u)}=y_{i-1}^{(u)}$ for all $u=1,\dots,k$ with probability 
at least $\e$. 
Thus (\ref{claa}) holds when $v=1$ and 
$\nep{-c} = 1-\e$. The case $v> 1$ is obtained by recursion.
Namely, the coupling $E_i^-$ can be further defined as follows. 
If we have a full matching $x_{i-1}^{(u)}=y_{i-1}^{(u)}$ for all
$u=1,\dots,k$ then we can match all the remaining particles. If we do
not have the full matching we sample particles labeled $i-2$ by
$\wt \nu$ with $j=i-1$ and $m=2\ell$ and a
suitable choice of the values $x_{i-1}^{(u)},y_{i-1}^{(u)}$ (the ones
that were sampled in the first step). Again, if we have a full
matching for particles labeled $i-2$ we can match all the remaining
particles labeled $i-3,\dots,i-m$. We repeat this procedure at later
steps. Since (\ref{lemmo1}) shows that there is a probability at least
$1-\nep{-c}$ to have a full matching at every step the bound
(\ref{claa}) follows.

\smallskip
To prove (\ref{lemmo1}) we are going to use the same argument as in
the proof of Lemma \ref{lemma_spiaccicato}, see the proof of
(\ref{iab}) in particular. 
We first sample the pairs $x_{r}^{(1)},y_{r}^{(1)}$,
$r=j-m,\ldots,j-1$ with an {\sl independent coupling} for the
corresponding marginals. Note that even if $x_{j}^{(1)}\neq
y_{j}^{(1)}$ the positions $x_{j-1}^{(1)},y_{j-1}^{(1)}$ have the same
ground state (i.e.\ minimal position) dictated by the ceiling $\xi$
and the (common) boundary conditions for the particles labeled
$1,\dots,j-m-1$ of polymer $1$. From the proof of Lemma
\ref{lemma_spiaccicato} we know that this implies that the matching
event $ x_{j-1}^{(1)}=y_{j-1}^{(1)}$ has probability at least $\d^2$
where $\d=\d(\a)>0$ is the probability that $x_{j-1}^{(1)}$ equals its
ground state position.  Next, we sample $x_{r}^{(2)},y_{r}^{(2)}$,
$r=j-m,\ldots,j-1$ with an independent coupling of the marginals
conditioned on the configuration of
$\{x_{r}^{(1)},y_{r}^{(1)}\}_{r=j-m,\ldots,j-1}$ which we extracted in
the previous step.
We claim that, conditionally on the occurrence of the matching
$x_{j-1}^{(1)}=y_{j-1}^{(1)}$, the matching event
$x_{j-1}^{(2)}=y_{j-1}^{(2)}$ 
has probability at least $\d^2$. 
Indeed, even if 
the ceilings  ``felt'' by the two copies of the positions
$x_{r}^{(2)},y_{r}^{(2)}$ of particles of polymer $2$ are in general distinct,
the ground state positions for $x_{j-1}^{(2)},y_{j-1}^{(2)}$ are dictated only
by the positions $x_{j-1}^{(1)},y_{j-1}^{(1)}$ which are now assumed
to coincide.  In particular the argument from the proof of Lemma \ref{lemma_spiaccicato}
again applies. This procedure can be repeated until the last polymer is
reached and in conclusion the independent coupling gives probability
at least $\e=\d^{2k}$ to the full matching event in (\ref{lemmo1}). This ends the proof of Lemma \ref{proel}.
\qed

\bigskip

In the proof of Lemma \ref{proel} (cf.\ (\ref{couple1})) we have used
a continuous time version of the so-called {\em path coupling}
argument, see \cite{BD} or e.g.\ \cite[Theorem 14.5]{cf:peres-book}
for the usual discrete time version.  For the sake of completeness we
give a proof in the next lemma.
\begin{Le}\la{path}
Suppose that (\ref{couple1}) holds for all pairs 
$\eta,\eta'$ such that $d_0(\eta,\eta')=1$.
Then 
$\bbE(d_\g(\eta(t),\eta'(t)))\leq k N(h,n)\nep{-t}$, for all $t\geq 0$ and for all initial data $\eta,\eta'$.
\end{Le}
\proof
We may define a graph having as vertices the elements of $E_{\xi,\si}$ by
declaring a pair $\eta,\eta'$ to be an edge whenever 
$d_0(\eta,\eta')= 1$.
For any $\eta,\eta'$, let $\O(\eta,\eta')$ denote the set of paths
connecting $\eta$ and $\eta'$, i.e.\ $\o\in\O(\eta,\eta')$ if
$\o=(\o_1,\dots,\o_r)$, $\o_1=\eta,\o_r=\eta'$ and
$d_0(\o_i,\o_{i+1})=1$, $i=1,\dots,r-1$. Then one checks that
\be\la{pathmetric} d_\g(\eta,\eta') =
\min_{\o\in\O(\eta,\eta')}\sum_{i=1}^{r-1}d_\g(\o_i,\o_{i+1})\,.
\end{equation}
To prove
(\ref{pathmetric}) observe that it suffices to exhibit one path which
achieves the equality since by construction it is clear that $d_\g
(\eta,\eta') \leq \sum_{i=1}^{r-1}d_\g(\o_i,\o_{i+1})$ for any
$\o\in\O(\eta,\eta')$.  Such a path can be informally defined as
follows. Consider the bottom paths $\eta^{(k)}, \eta'^{(k)}$ and the
positions $x$ where $\eta^{(k)}_x>\eta'^{(k)}_x$. Then move one by one
the particles of $\eta^{(k)}$ in this region until we have
$\eta^{(k)}_x\leq \eta'^{(k)}_x$ everywhere. Then consider the paths
$\eta^{(k-1)}, \eta'^{(k-1)}$ and the positions $x$ where
$\eta^{(k-1)}_x>\eta'^{(k-1)}_x$. As before, move one by one the
particles of $\eta^{(k-1)}$ in this region until we have
$\eta^{(k-1)}_x\leq \eta'^{(k-1)}_x$ everywhere. We repeat this
procedure until we reach the top paths.  At this point we have reached
a configuration $\wt\eta$ such that $\wt \eta \leq \eta$ (everywhere).
Next we start from the top paths $\eta^{(1)}, \eta'^{(1)}$ and
consider the region where $\wt \eta^{(1)}_x<\eta'^{(1)}_x$. We can
move one by one the particles in this region until we have $\wt
\eta^{(1)}=\eta'^{(1)}$ (everywhere). We repeat with the paths $\wt
\eta^{(2)},\eta'^{(2)}$ and so on until we reach the bottom polymers
labeled $k$. This construction produces a path which realizes the
minimum (\ref{pathmetric}) since we never used more than the strictly
necessary moves.

From the triangle inequality, 
for each pair of initial conditions $\eta,\eta'$ we have
$$
\bbE[d_\g(\eta(t),\eta'(t))]\leq \sum_{i=1}^{r-1}
\bbE[d_\g(\zeta_i(t),\zeta_{i+1}(t)]\,,
$$
where we call $\zeta_1,\dots,\zeta_r$ the minimizing path in (\ref{pathmetric}) and $\zeta_1(t),\dots,\zeta_r(t)$ the corresponding trajectory. 
In particular, subtracting $d_\g(\eta,\eta')=\sum_{i=1}^{r-1}d_\g(\zeta_i,\zeta_{i+1})$, 
dividing by $t$ and letting $t\downarrow 0$ we obtain
$$
[\wt\cG d_\g](\eta,\eta')\leq \sum_{i=1}^{r-1} [\wt \cG d_\g](\zeta_i,\zeta_{i+1})\,.
$$
Since each term in the r.h.s.\ above is of the form 
(\ref{couple1}) with $d_0(\zeta_i,\zeta_{i+1})=1$
we have that the assumptions imply
$$
[\wt \cG d_\g](\eta,\eta')\leq - d_\g(\eta,\eta')\,,
$$ 
for arbitrary initial conditions. 
Therefore $\varphi(t):=\bbE[d_\g(\eta(t),\eta'(t))]$ satisfies 
$\frac{d}{dt}\varphi(t) \leq - \varphi(t)$, which implies the claim
since $\varphi(0)\le \max_{\eta,\eta'}d_\gamma(\eta,\eta')\le k N(h,n)$.
\qed

\subsection{A dynamics with full-polymer moves}\la{fulll}
The next ingredient which enters the proof of Theorem \ref{teo_gap} is
a dynamics where each move consists in updating $(2s+1)$ whole
polymers, $s\in\N$.  As usual, we let $\xi,\sigma\in\Omega^h_{1,n}$
(with $\sigma\leq \xi$) and $\mu^{\xi,\sigma}(\cdot)$ denotes the law
on $E_{\xi,\sigma}$ for $k$ polymers with floor and ceiling
$\sigma,\xi$, as defined in \eqref{eq1}.  To each $1\leq j\leq k$ is
assigned an independent Poisson clock of mean $1$. When the clock
labeled $j$ rings, we update the polymers $\eta^{(u)}$ with index $
\max(j-s,1)\leq u\leq \min(j+s,k)$, sampling the new configuration
according to the law
\begin{eqnarray}
\label{eq:hatnu}
  \hat \nu_{j,s}(\cdot):=  
\mu^{\xi,\sigma}\left(\,\cdot\,\tc\eta^{(\max(j-s-1,0))},
    \eta^{(\min(j+s+1,k+1))}\right),
\end{eqnarray}
with the convention that $\eta^{(0)}:=\xi$ and $\eta^{(k+1)}:=\sigma$.
Call $\mathcal M$ the generator of this dynamics.

\begin{Proposition}\la{profull}
  For every $\alpha>0$ there exists $s:=s(\alpha)\in\N$ such that,
uniformly in $n,k$, $|h|\leq n$  and on $\sigma\leq \xi$, one has
\begin{eqnarray}
\label{eq:gapM}
  \gap(\mathcal M)\geq 1.
\end{eqnarray}
\end{Proposition}
\proof
The general structure of the proof is similar to that of Proposition
\ref{pros}, but the coupling argument is rather different. 
Given $\rho>0$, we define the distance 
function $D_\rho(\cdot,\cdot)$ by
setting 
for every $\eta,\eta'$,
\begin{eqnarray}
\label{eq:Drho}
  D_\rho(\eta,\eta'):=\sum_{j=1}^k \nep{-j\rho}\sum_{x=1}^{n-1}
\frac{|\eta^{(j)}_x-\eta'^{(j)}_x|}2.
\end{eqnarray}
In analogy with \eqref{pathmetric}, one checks
that 
\begin{eqnarray}
  D_\rho(\eta,\eta')=\min_{\o\in\Omega(\eta,\eta')}\sum_{i=1}^{r-1}
D_\rho(\o_i,\o_{i+1}),
\end{eqnarray}
where in this case one requires that $D_0(\o_i,\o_{i+1})=1$ for $i<r$.

Given two initial conditions $(\eta,\eta')$, let $(\eta(t),\eta'(t))$
be the corresponding evolutions. As in the proof of Proposition
\ref{pros}, cf.\ Lemma \ref{path}, to prove \eqref{eq:gapM}, it is then
sufficient to prove that for every $\alpha>0$ there exists a choice of
$s\in\N$, $\rho>0$ and a coupling $\bbP$ of $(\eta(t),\eta'(t))$ such
that
 \begin{eqnarray}
\label{auxgapM}
   \left[\wt{\mathcal M} D_\rho\right](\eta,\eta'):=
\left.\frac{d}{dt}\bbE(D_\rho(\eta(t),\eta'(t)))\right |_{t=0}\leq
-D_\rho(\eta,\eta')
\end{eqnarray}
whenever $D_0(\eta,\eta')=1$.

Let $(\bar\eta,\bar\eta')$ satisfy the
 latter condition, with the single discrepancy consisting in
 $\bar\eta^{(i)}_x=\bar\eta'^{(i)}_x+2$.  In analogy with
 \eqref{upto}, we write the generator of the coupled dynamics as
\begin{eqnarray}
  \wt {\mathcal M}:=\sum_{j=1}^k\left(
\hat E_{j,s}-1\right),
\end{eqnarray}
where $\hat E_{j,s}$ is a coupling (to be specified) of $\hat
\nu_{j,s}$ for the two configurations. For all $j$ such that
$i-s\leq j\leq i+s$ we can choose the coupling such that 
$[\hat E_{j,s}D_\rho](\bar\eta,\bar\eta')=0$. Moreover,
 if $j\in\{1,\ldots,k\}\setminus\{i-s-1,\ldots,i+s+1\}$
we can choose the coupling such that 
$[\hat E_{j,s}D_\rho](\bar\eta,\bar\eta')=D_\rho(\bar\eta,\bar\eta')$.
 One has therefore
\begin{align}\label{dup3}
&   [\wt {\mathcal M}D_\rho](\bar\eta,\bar\eta')\\
&\quad\leq 
-s D_\rho(\bar\eta,\bar\eta')+
  [(\hat E_{i+s+1,s}-1)D_\rho](\bar\eta,\bar\eta')+
  [(\hat E_{i-s-1,s}-1)D_\rho](\bar\eta,\bar\eta').
  \nonumber
  \end{align}
It is clear that the last two terms may be non-negative, and that they vanish 
if $i\geq k-s$ or $i\leq s+1$, respectively. 

Let us analyze first the easier case of $\hat E_{i+s+1,s}$ in which case, 
it is worth recalling, we are updating polymers labeled $i+1,\ldots,
i+2s+1$. Since
$\bar\eta^{(i)}\geq\bar\eta'^{(i)}$ while
$\bar\eta^{(i+2s+2)}=\bar\eta'^{(i+2s+2)}$, by monotonicity there exists a
coupling $\hat E_{i+s+1,s}$ such that one has $\eta^{(j)}\geq
\eta'^{(j)}$ for every $i<j<i+2s+2$. Moreover, since $\bar\eta^{(i)}$
and $\bar\eta'^{(i)}$ differ only at $x$, we can choose 
$\hat E_{i+s+1,s}$ such that $\eta^{(j)}_y=\eta'^{(j)}_y$ for every
$i<j<i+2s+2$ and $y$ outside the interval $\{a_-,\ldots,a_+\}$, where
$a_+:=\inf\{y>x: \eta'^{(i+2s+1)}_y=\bar\eta^{(i)}_y\}$ and
$a_-:=\sup\{y<x:\eta^{(i+2s+1)}_y=\bar\eta^{(i)}_y\}$.
As a consequence, going back to the definition of $D_\rho(\cdot,\cdot)$,
\begin{eqnarray}
  [\hat E_{i+s+1,s}D_\rho](\bar \eta,\bar\eta')\leq
  \sum_{j=i+1}^{\infty}e^{-j\rho}\hat\nu_{i+s+1,s}\left((a_+-a_-)^2\right).
\end{eqnarray}
Thanks to  \eqref{eq:spiacc2} there exists $c_3(\alpha)>0$ such 
that for every $u>0$
\begin{eqnarray}
\label{eq:a+a-}
\hat\nu_{i+s+1,s}(a_+-a_-=u)\leq e^{-c_3(\alpha)u}.  
\end{eqnarray}
From this one deduces immediately that there exists $c_4(\alpha,\rho)<\infty$
such that
\begin{eqnarray}
\label{eq:dup}
    [\hat E_{i+s+1,s}D_\rho](\bar \eta,\bar\eta')\leq c_4(\alpha,\rho)D_\rho
(\bar\eta,\bar\eta').
\end{eqnarray}

\medskip

Finally we deal with $\hat E_{i-s-1,s}$. We have from \eqref{eq:Drho}
\begin{eqnarray}
\label{eq:Emeno}
[  \hat E_{i-s-1,s}D_\rho](\bar \eta,\bar \eta')=\frac12
\sum_{j=i-2s-1}^{i-1}e^{-\rho j}
\hat E_{i-s-1,s}\left(
\sum_{y=1}^{n-1}|\eta^{(j)}_y-\eta'^{(j)}_y|
\right).
\end{eqnarray}
Again, we can choose the coupling such that $\eta'^{(j)}\leq
\eta^{(j)}$ and $\eta^{(j)}_y=\eta'^{(j)}_y$ for
$y\notin\{b_-,\ldots,b_+\}$, where
$b_-=\sup\{y<x:\eta'^{(i-1)}_y=\bar \eta^{(i-2s-2)}_y\}$, and 
similarly for $b_+$. In analogy with \eqref{eq:a+a-} one has 
\begin{eqnarray}
\label{eq:b+b-}
\hat\nu_{i-s-1,s}(b_+-b_-=u)\leq e^{-c_3(\alpha)u}.  
\end{eqnarray}
  Then,
\begin{eqnarray}
\label{sumetaaa}
\sum_{y=1}^{n-1}|\eta^{(j)}_y-\eta'^{(j)}_y|\leq 
|b_+-b_-|\times|\{b_-<y<b_+:\;\eta'^{(j)}_y
\ne \bar\eta^{(i-2s-2)}_y\}|.
\end{eqnarray}
Using \eqref{sumetaaa}, \eqref{eq:spiacc1} and \eqref{eq:b+b-},
\begin{align*}
\hat   E_{i-s-1,s}\left(
    \sum_{y=1}^{n-1}|\eta^{(j)}_y-\eta'^{(j)}_y|
  \right)&=\hat
  E_{i-s-1,s}\left[\left.\hat E_{i-s-1,s}
\left(\sum_{y=1}^{n-1}|\eta^{(j)}_y-\eta'^{(j)}_y|
      \right)\right|b_-,b_+\right]
\\
&\leq \hat\nu_{i-s-1,s}\left[(b_+-b_-)^2\right]e^{-c_4(\alpha)(i-j)}
\leq c_5(\alpha)e^{-c_4(\alpha)(i-j)}.
\end{align*}
Therefore, going back to \eqref{eq:Emeno},
\begin{eqnarray}
\label{dup2}
[  \hat E_{i-s-1,s}D_\rho](\bar \eta,\bar \eta')\leq
c_5 (\alpha)D_\rho(\bar\eta,\bar\eta')\sum_{r=1}^\infty e^{-(c_4(\alpha)-\rho)r}=
c_6(\alpha)  D_{\rho(\alpha)}(\bar\eta,\bar\eta'),
\end{eqnarray}
where we chose $\rho:=\rho(\alpha)=c_4(\alpha)/2$. From \eqref{dup3},
\eqref{eq:dup} and \eqref{dup2} one concludes that
\begin{eqnarray}
 [\wt {\mathcal M}D_{\rho(\alpha)}](\bar\eta,\bar\eta')\leq   
-(s-c_7(\alpha)) D_{\rho(\alpha)}(\bar\eta,\bar\eta').
\end{eqnarray}
At this point, it is sufficient to choose $s:=s(\alpha):=\lceil c_7(\alpha)
\rceil+1$ to get \eqref{auxgapM}.
\qed

\subsection{An estimate for $k$ polymers with $r$ particles}\la{bounded}
The last ingredient we need for the proof of Theorem \ref{teo_gap} is
a rough estimate on the spectral gap for a system with $k$ polymers,
each with $r$ particles; it is important that this bound is
independent of the lengths $n_1,\dots,n_k$ of each polymer.  Consider
a configuration $\eta\in E_{\xi,\si}\subset\O^h_{k,n}$ with $k$
polymers, each with $N$ particles. Let ${\bf x}={\bf x}(\eta)$ denote
the corresponding particle configuration. Fix $0\leq j<m\leq N+1$ and
consider the probability measures $\nu^{{\bf x}}_{j,m}$ defined in
(\ref{partcond2}).  If we freeze all particles labeled $i\leq j$ and
$i\geq m$ we can perform the local-update dynamics defined in
(\ref{diro}) for the $r:=m-j-1$ particles labeled $j+1,\dots,m-1$.
This process is clearly reversible with respect to $\nu^{{\bf
    x}}_{j,m}$.  Its Dirichlet form is given by
$$
\cE^{\bf x}_{j,m}(f)= \sum_{u=1}^k\sum_{x=x_j^{(u)}+1}^{x_m^{(u)}-1}
\nu^{{\bf x}}_{j,m}\left[\var_{u,x}(f) \right]\,,
$$
where $\var_{u,x}(f)$ has the same meaning as in (\ref{diro}).  As
usual, below we use the notation $\var_{\nu}(f)=\nu(f^2)-\nu(f)^2$ for
any probability measure $\nu$.
\begin{Le}\la{crude}
  For every $\a>0$, $k\in\bbN$ and $r\in\bbN$, there exists
  $c=c(\a,k,r)$ such that for any pair $j,m$ with $m-j-1=r$, for any
  choice of all other parameters and for any $f:\O^h_{k,n}\to\bbR$
$$
\cE^{\bf x}_{j,m}(f)\geq c\,\var_{\nu^{{\bf x}}_{j,m}}(f)\,.
$$
\end{Le}
\proof The only delicate point here is that the length of the portion
of paths where the $r$ particles live is arbitrarily long and we need
an estimate which does not depend on that.  We introduce a further
family of measures as follows.  Let $\nu^{{\bf x},u}_{j,m}$ denote the
law of $\eta^{(u)}$ according to the measure $\nu^{{\bf x}}_{j,m}$
conditioned on the value of the paths $\eta^{(i)},\,i\neq u$.  If we
let $$\cE^{{\bf x},u}_{j,m}(f)= \sum_{x=x_j^{(u)}+1}^{x_m^{(u)}-1}
\nu^{{\bf x},u}_{j,m}[\var_{u,x}(f)]\,$$ denote the corresponding
Dirichlet form we know from the $k=1$ version of Theorem \ref{teo_gap}
(see the remark following the statement of Theorem \ref{teo_gap}) that
\be\la{gap10} \cE^{{\bf x},u}_{j,m}(f)\geq \d(\a)\,\var_{\nu^{{\bf
      x},u}_{j,m}}(f)\,,
\end{equation}
for some constant $\d(\a)$ depending only on $\a$.  Taking expectation
w.r.t.\ $\nu^{{\bf x}}_{j,m}$, using $\nu^{{\bf x}}_{j,m}[\nu^{{\bf
    x},u}_{j,m}(g)]=\nu^{{\bf x}}_{j,m}(g)$ for any function $g$, and
summing over $u$ in (\ref{gap10}), by definition of $\cE^{{\bf
    x}}_{j,m}(f)$ we have
$$
\cE^{{\bf x}}_{j,m}(f) \geq  \d(\a)
\sum_{u=1}^k\nu^{{\bf x}}_{j,m}\left[\var_{\nu^{{\bf x},u}_{j,m}}(f)
\right]\,.
$$
Therefore it remains to prove that for some $c=c(\a,r,k)>0$ one has \
\be\la{gap11} \sum_{u=1}^k \nu^{{\bf x}}_{j,m}\left[\var_{\nu^{{\bf
        x},u}_{j,m}}(f) \right]\geq c\,\var_{\nu^{{\bf
      x}}_{j,m}}(f)\,,
\end{equation}
for all functions $f$. To prove this estimate we observe that the
left-hand side of (\ref{gap11}) coincides with the Dirichlet form of
the Markov chain described as follows. Attach independent rate 1 Poisson
clocks to the labels $u=1,\dots,k$. When a label $u$ rings update the
whole path between $x_j^{(u)}+1$ and $x_m^{(u)}-1$ according to the
distribution $\nu^{{\bf x},u}_{j,m}$ (that is, freeze all other
polymers and update polymer $u$ with a sample from $\nu^{{\bf
    x},u}_{j,m}$).

Thus, the following rough coupling argument will suffice for the proof
of (\ref{gap11}). Namely, consider the Markov chain started in the
minimal configuration (i.e.\ each of the $k$ polymers starts in the
minimal path compatible with the particles $x_j^{(u)}$ and
$x_m^{(u)}$). Let $E_t$ denote the event that up to time $t$ the
Markov chain has never visited the maximal allowed configuration.  It
is not hard to prove a bound of the form $\bbP(E_t)\leq
c^{-1}\nep{-c\,t}$ for some $c=c(\a,r,k)>0$.  This, in turn, implies
the desired spectral gap estimate (\ref{gap11}) using monotonicity and
the bound (\ref{gapcoup}).

To prove $\bbP(E_t)\leq c^{-1}\nep{-c\,t}$ we may reason as follows.
Consider the event $F_s$ that in the time interval $[s,s+1)$ the
following sequence of $k$ successive updatings appears: for
$i=1,\dots,k$ the $i$-th update is for polymer $u=i$ and the
configuration sampled from $\nu^{{\bf x},u}_{j,m}$ is the maximal
allowed path for the $i$-th polymer given the current configuration.
Since there are $r$ particles, an application of Lemma
\ref{lemma_spiaccicato} shows that $\bbP(F_s)\geq p$ where
$p=p(\a,r,k)>0$ depends 
neither the lengths
$n_u=x_m^{(u)}-x_j^{(u)}$ of the polymers, nor on the configuration
of $\eta$ at time $s$. By construction, on the
event $F_s$ we know that the Markov chain visited the maximal allowed
configuration at least once in the time interval $[s,s+1)$. Therefore
the event $E_t$ implies that none of the events $F_s$ occurred for
$s=0,\dots,\lfloor t-1\rfloor$.  This implies $\bbP(E_t)\leq (1-p)^{t-1}\leq
(1-p)^{-1}\nep{-p\,t}$.  The proof of Lemma \ref{crude} is complete.
\qed

\subsection{Putting everything together: proof of Theorem \ref{teo_gap}}\la{together}

Once Proposition \ref{pros}, Proposition \ref{profull} and Lemma
\ref{crude} are established, the proof of Theorem \ref{teo_gap} is
obtained through a chain of comparison inequalities.  Indeed,
Proposition \ref{profull} may be restated as 
\begin{equation}
\la{gap13} \var(f)\leq
\sum_{j=1}^k\mu\left[\var_{\hat\nu_{j,s}}(f)\right]\,,
\end{equation}
where $\var(f)$ is the variance w.r.t.\ $\mu:=\mu^{\xi,\si}$ and
$\hat\nu_{j,s}$ is the conditional probability measure defined in
\eqref{eq:hatnu}.
Here $s=s(\a)$ is a fixed parameter ($2s+1$ is the number of polymers
to be updated at each step in the process with generator $\cM$
appearing in Proposition \ref{profull}: they are the polymers labeled $u$, 
with
$\max(j-s-1,1)\le u\le\min(j+s+1,k)$).  For each $j$ in
(\ref{gap13}) we apply the estimate of Proposition \ref{pros} with the
number of polymers equal to $k=s$.  This yields 
\begin{equation}
  \la{gap14}
  \var_{\hat\nu_{j,s}}(f)\leq \sum_{i=1}^N \hat\nu_{j,s}\left[
    \var_{\r^{{\bf
    x}}_{i,\ell}}(f)\right] \,,
\end{equation}
where $\r^{{\bf x}}_{i,\ell}$ is the conditional probability measure
appearing in (\ref{dirrho}) and $\ell=\ell(\a)$ is a fixed parameter
($2\ell+1$ is the size of the block of particles to be updated at each
step in the process with generator $\cG$ appearing in Proposition
\ref{pros}).  Finally, recalling the definition $\r^{{\bf
    x}}_{i,\ell}=\nu^{{\bf x}}_{i-\ell-1,i+\ell+1}$ (cf.\
(\ref{partcond2})) and applying the bound in Lemma \ref{crude} (with
$r=2\ell+1$ and $k=s$) we know that for every $i$ in (\ref{gap14}) we
have the estimate for all ${\bf x}={\bf x}(\eta)$: 
\begin{equation}
  \la{gap15}
  \var_{\r^{\bf x}_{i,\ell}}(f)\leq c^{-1}\sum_{u=j-s-1}^{j+s+1}
  \sum_{x=x_{i-\ell-1}^{(u)}+1}^{x_{i+\ell+1}^{(u)}-1} \r^{\bf
    x}_{i,\ell}\left[\var_{u,x}(f) \right]\,.
\end{equation}
Note that the constant $c=c(\a,s,\ell)$  depends only on $\a$ here.
Also, note that since there are at most $2\ell+1$ blocks of particles
covering a given position $x$ in the path we obtain 
\begin{equation}
  \la{gap25}
  \sum_{i=1}^N \hat\nu_{j,s}\left[ \var_{\r^{{\bf x}}_{i,\ell}}(f)\right] \leq
  c^{-1}(2\ell+1) \sum_{u=j-s-1}^{j+s+1} \sum_{x=1}^{n-1}
  \hat\nu_{j,s} \left[\var_{u,x}(f) \right]\,,
\end{equation}
where we use the property of conditional expectation $ \hat\nu_{j,s}\left[ 
\r^{{\bf x}}_{i,\ell}(g)\right] = \hat\nu_{j,s}(g)$ valid for any function $g$.

Putting everything together and using
$\mu\left[\hat\nu_{j,s}(g)\right] =\mu(g)$ for any $g$, from
(\ref{gap13})-(\ref{gap25}) we see that
\begin{align*}\la{gap16}
\var(f)&\leq c^{-1}(2\ell+1) 
\sum_{j=1}^k\sum_{u=j-s-1}^{j+s+1}
\sum_{x=1}^{n-1}
\mu\left[\var_{u,x}(f)
\right]\\
&\leq c^{-1}(2\ell+1)(2s+1) 
\sum_{i=1}^k\sum_{x=1}^{n-1}
\mu\left[\var_{i,x}(f)
\right] \\&=c^{-1}(2\ell+1)(2s+1)\cE(f) \,.
\end{align*}
This shows that $\gap(\a)\geq c'(\a)$, with $c'(\a)=c\,(2\ell+1)^{-1}(2s+1)^{-1}$, which 
completes the proof of Theorem \ref{teo_gap}.
\qed

\section{Proof of Theorem \ref{th:tmix1}}\la{mixing}
For clarity of exposition we give the proof first of all in the
simpler case where the ceiling $\xi$ is the maximal configuration, $\wedge$, in
$\Omega^h_{1,n}$. Later, in
Section \ref{sec:xisigma} we sketch the modifications which are needed
in the general situation.


\subsection{The case of maximal ceiling}
\label{sec:maxceil} Fix some positive $T>0$.
By monotonicity and the definition \eqref{eq:Tmix} of mixing time, it
is clear that a sufficient condition for $\tmix\leq T$ is that the
first (random) time when the dynamics started from the minimal
configuration, $\{\eta^{(j)}=\si,j=1,\ldots,k\}$, hits the
maximal configuration, $\{\eta^{(j)}=\xi, j=1,\ldots,k\}$, is smaller
than $T$ with probability at least $1-1/(2e)$.  Again by monotonicity,
it is easy to convince oneself that this random time does not decrease
if one replaces $h$ with zero, $n$ with $2M:=2\max(n,k)$, $k$ with $M$
and $\si$ with $\vee$.

\medskip

Therefore, to prove Theorem \ref{th:tmix1} (in the case of maximal
ceiling) it is sufficient to prove the following.  Let $h=0$, $n=2M$, $k=M$,
$\xi=\wedge$, $\si=\vee$, cf. \eqref{eq:wedge}, \eqref{eq:vee}; start  
the dynamics from the minimal configuration $\eta_-:=\{\eta^{(j)}=\vee
$ for all $
j=1,\ldots,M\}$, and call $t(M)$ the first time when the 
maximal configuration $\eta^+:=\{\eta^{(j)}=\wedge$ for all $
j=1,\ldots,M\}$ is reached. Then,
\begin{Theorem}
\label{th:ineqt}
There exists $C(\alpha)>0$ such that
with probability larger than $1-1/(2e)$
\begin{eqnarray}
\label{eq:ineqt}
t(M)\leq C(\alpha)M\,(\log M)^6.   
\end{eqnarray}
\end{Theorem}
Theorem \ref{th:ineqt} will be proven just after the statement of
Theorem \ref{th:souffle'} below. Before stating Theorem
\ref{th:souffle'}, we need to introduce in Sections  \ref{sec:polyset}
 and \ref{sec:proofof} a few
auxiliary definitions.
\medskip

{\sl A notational convention} When in the following we say that an
event occurs ``with large probability'' (w.l.p.), we mean to say that
the probability of the complementary event goes to zero for
$M\to\infty$ faster than any inverse power of $M$.  Since, as will be
clear, we have to exclude only polynomially many (in $M$) events which w.l.p.\ do not occur,
by the union bound we have that the occurrence of at least one of
these events still goes to zero for $M\to\infty$.  For simplicity of
exposition, and when there is no risk of confusion, we will often just
pretend that an event which occurs w.l.p., occurs deterministically.

\subsubsection{From polymer configurations to subsets of the cube}
\label{sec:polyset}
\begin{figure}[h]
\centerline{
\psfrag{x}{$y$}
\psfrag{y}{$x$}
\psfrag{0}{$0$}
\psfrag{m}{$M$}
\psfrag{2M}{$2M$}
\psfig{file=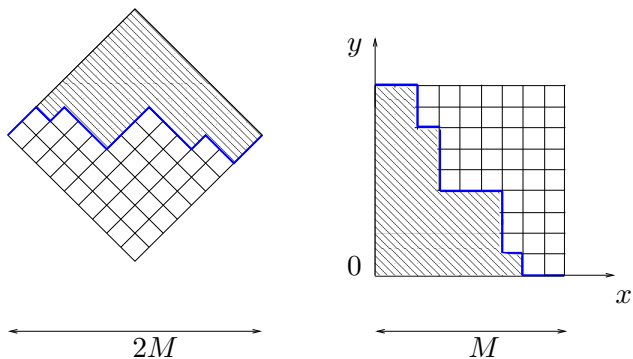, width=3.2in}}
\caption{A graphical construction of the set $s(\eta)$. Given $\eta\in
  \Omega^0_{M,2M}$, in order to obtain $s(\eta)$ do the following, for
  $r_3=0,\ldots,(M-1)$: a) draw the configuration of $\eta^{(M-r_3)}$
  (thick line in the left drawing); b) rotate the picture by $3\pi/4$
  anti-clockwise and shrink it by a scale factor $1/\sqrt 2$. The
  shaded region coincides with the horizontal section of $s(\eta)$ at
  height $r_3<h\leq r_3+1$. It is obvious from this construction that
  the subset $s(\eta)$ thus obtained is a monotone subset of the cube
  $\mathcal C_M$.  }
\label{fig:seta}
\end{figure}

The proof of Theorem \ref{th:ineqt} becomes more intuitive if one
interprets a configuration $\eta\in\Omega^0_{M,2M}$ as a subset $s$ of
the cube $\mathcal C_M:=[0,M]^3\subset \bbR^3$. This mapping is just another way to see the 
mapping,  introduced
in Section \ref{sec:fromnon},  between polymer configurations $\eta$ and boxed plane partitions. To define precisely
this correspondence, divide first of all $\mathcal C_M$ into $M^3$
elementary cubes of unit side, which we label with the (integer)
coordinates $r=(r_1,r_2, r_3)$ of their point of smallest $L^1$ norm
(observe that $0\leq r_i<M)$.  Then, to a given
$\eta\in\Omega^0_{M,2M}$ we associate $s=s(\eta)\subset \mathcal C_M$,
a union of elementary cubes, by establishing that the elementary
cube labeled $r$ belongs to $s$ if and only if (cf. Figure \ref{fig:seta})
\begin{eqnarray}
  \eta^{(M-r_3)}_{M-r_1+r_2}<M-r_1-r_2.
\end{eqnarray}

It is easy to check that:
\begin{itemize}
\item if $\eta\leq \eta'$ then $s'\subset s$ (note that the inequality is
reversed!)
\item the maximal configuration $\eta^+$ defined above corresponds to
  the empty subset of the cube: $s^-=\emptyset$. Conversely, the
  minimal polymer configuration $\eta^-$ corresponds to the maximal
  subset $s^+=\mathcal C_M$, i.e., the full cube. We will sometimes
  refer to $s^-$ as to the {\sl ground state}, for obvious reasons.
\item for every $\eta\in\Omega^0_{M,2M}$, $s(\eta)$ is a {\sl monotone
    subset} of $\mathcal C_M$, i.e., if the elementary cube labeled $r$
  belongs to $s$, then so do also all the elementary cubes $r'$ such that
  $r'_1\leq r_1,r'_2\leq r_2, r'_3\leq r_3$.
\end{itemize}

The equilibrium measure is described in terms
of $s$ by: 
\begin{eqnarray}
  \mu^{\wedge,\vee}(s)=\frac{e^{-2\alpha|s|}}{\hat Z_{\wedge,\vee}},
\end{eqnarray}
where $|s|$ denotes the number of elementary cubes contained in $s$,
i.e., its volume.  

The dynamics of Section \ref{sec:results} can also be explicitly described in terms
of $s$. Here, let us simply remark that the elementary moves
of the Markov Chain consist in adding or removing a single elementary cube,
with the constraint that $s$ remains a monotone subset of $\mathcal C_M$ after the update.
Observe also that, under the $\eta\leftrightarrow s$ correspondence,
the upward drift felt by the polymers $\eta^{(j)}$ during the dynamics translates into the
fact that the upper boundary of $s$ feels a drift in the direction
$(-1,-1,-1)$.

\bigskip

\subsubsection {An auxiliary dynamics} 
\label{sec:proofof}
To avoid a plethora
of $\lfloor \cdot \rfloor$, we assume that $(\log M)^2\in 2\N$ and
that $K:=(M/(\log M)^2)\in\N$.  Divide $\mathcal C_M$ into sub-cubes
$B_v$ (called {\sl blocks} from now on) of side $(\log M)^2$, indexed
by $v=(v_1,v_2,v_3)$ with $0\leq v_i<K$, and such that the point of
$B_v$ with minimal $L^1$ norm is $(\log M)^2 v$ (of course, $v_i$ are
integers).  Given $v$, we will also define $B_v^+,B_v^-$ to be the
half-blocks obtained cutting $B_v$ horizontally into two equal parts
($B_v^-$ will denote the bottom one).  Call $s_t$ the configuration at
time $t$, which starts from the completely full configuration $s^+$ at
time $t=0$.

The idea behind Theorem \ref{th:ineqt} is to consider an easier dynamics $\hat s_t$ such that
$s_t\subset\hat s_t$ almost surely and, calling $\hat t(M)$ the first
time $\hat s_t$ reaches the empty configuration $s^-$, to show that
$\hat t(M)$ satisfies \eqref{eq:ineqt} with probability at least
$1-1/(2e)$. The statement of  Theorem \ref{th:ineqt}  then follows
immediately by monotonicity, since $t(M)\le \hat t(M)$.

Let 
\begin{eqnarray}
  \label{eq:tau}
\tau:=\tau(M):=C(\alpha)(\log M)^{8}/7,
\end{eqnarray}
where $C(\alpha)$ is
the same as in \eqref{eq:ineqt}. We define now two {\sl deterministic}
sets $S^\pm_t\subset \mathcal C_M$ which (roughly speaking) coincide with 
$s^+$ at $t=0$, are empty after time $C(\alpha) M(\log
M)^6$, and such that $\hat s_t$ 
satisfies $S^-_t\subset \hat s_t\subset S^+_t$, w.l.p.\ and for all
$t\leq M^2$
(the latter property is non-trivial and it is the content of Theorem
\ref{th:souffle'} below).

\begin{figure}[!h]
\centerline{
\psfrag{x}{$x$}
\psfrag{y}{$y$}
\psfrag{z}{$z$}
\psfrag{s}{$\star$}
\psfig{file=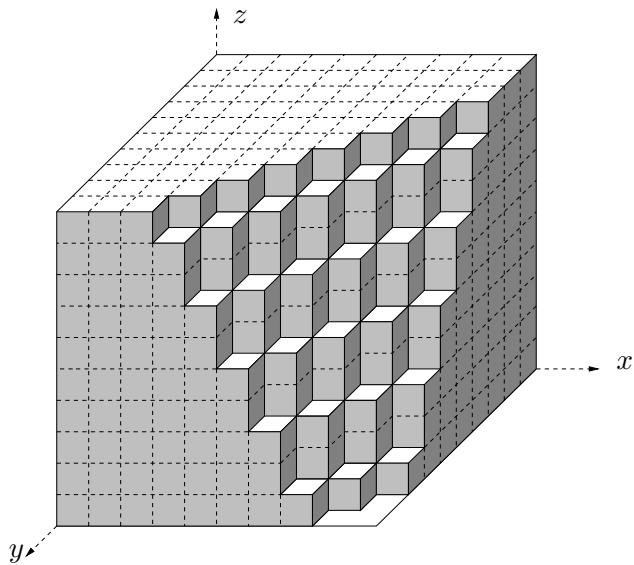,height=2.9in, width=3.2in}}
\caption{The set $S_t^-$ for $14\tau<t\leq 15\tau$. Small cubes denote
  blocks $B_v$ of side $(\log M)^2$, and there are $K=M/(\log M)^2(=10)$ of
  them along each side. At time intervals of $2\tau$, a new diagonal
  set of columns with $v_1+v_2=const$ starts to move downwards: from
  then on, it moves one block down each time interval $\tau$.  Once a
  column is empty, it stays empty forever. Roughly speaking, 
the set  $S_t^-$ contains all the blocks which are below a plane perpendicular
to the vector $(2,2,1)$ and which moves at constant speed, of order 
$(\log M)^{-6}$, in the direction 
$(-2,-2,-1)$.
The set $S^+_t$ for
  $i\tau\leq t<(i+1)\tau$ can be obtained simply by taking $S_t^-$ for
  some $(i-1)\tau<t\leq i\tau$ and adding a half-blocks $B_v^-$ on top
  of each incomplete but not empty column, and also to each empty
  column which is adjacent to a non-empty one ($S_t^+$ is not drawn in
  the picture). }
\label{fig:S-}
\end{figure}

$\{S^-_t\}_{ t\geq0}$, is defined as follows (see also Figure
\ref{fig:S-} for a graphical definition):
\begin{itemize}
\item $S^-_0=\mathcal C_M$
\item if $i\in \N$ and $(i-1)\tau<t\leq i\tau$, then $S_t^-$ contains all and 
only the blocks $B_v$ such that $v$ satisfies
\begin{eqnarray}
v_3+2(v_1+v_2)\leq 5(K-1)-i.
\end{eqnarray}
\end{itemize}
On the other hand, letting for ease of notation
$$
V_t^\pm:=\mathcal C_M\setminus S_t^\pm,
$$
$\{S^+_t\}_{ t\geq 0}$ is defined as follows (see also the caption of
Figure \ref{fig:S-}):
\begin{itemize}
\item $S^-_t\subset S^+_t\subset \mathcal C_M$
\item If $B_v\in V_{\tau\lfloor t/\tau\rfloor}^-$ and $B_{v-e_i}\in
  V_{\tau\lfloor t/\tau\rfloor}^-$ for at least one choice of $i=1,2,3$, then
  $B_v\in V^+_t$ ($e_i$ are the canonical base vectors of $\Z^3$, and 
$\lfloor x\rfloor:=\max\{n\in\bbZ:n\leq x\}$)
\item  If $B_v\in V_{\tau\lfloor t/\tau\rfloor}^-$ but there is no $i=1,2,3$ such that 
$B_{v-e_i}\in V_{\tau\lfloor t/\tau\rfloor}^-$, then $B_v^+\in  V^+_t$ but 
$B_v^-\in  S^+_t$.
\end{itemize}

The following properties of $S^\pm_t$ are immediately checked:
\begin{itemize}
\item $S_t^-$ is the union of blocks $B_v$, while $S_t^+$ is the union of
blocks {\sl and} of half-blocks $B_v^-$
\item  $S^-_\cdot$ is left-continuous, non-increasing, constant
in the time intervals $((i-1)\tau,i\tau]$ 
\item $S^+_\cdot$ is right-continuous, non-increasing and constant on
  intervals $[(i-1)\tau,i\tau)$
\item $S^\pm_t$ are monotone subsets of $\mathcal C_M$ for every $t$.
\item For  $t\geq (6/7)C(\alpha) M(\log M)^6=6K\tau$, $S^-_t$ is empty and
$S_t^+$ contains only 
$B^-_{(0,0,0)}$ ($K:=M/(\log M)^2$ and $\tau=C(\alpha)(\log M)^8/7$ were defined at the beginning of this section).
\end{itemize}

As an example, if $0<t\leq \tau$ then $S^-_t=\mathcal C_M\setminus
B_{(K-1,K-1,K-1)}$, while if $2\tau<t\leq 3\tau$
then
$$S^-_{t}=\mathcal C_M\setminus(
B_{(K-2,K-1,K-1)}\cup B_{(K-1,K-2,K-1)}\cup_{i=1}^3 B_{(K-1,K-1,K-i)}).$$

We define the auxiliary dynamics $\hat s_t$ by establishing
that it  has the same law as
$s_t$ {\sl conditioned on the event that $S^-_u\subset s_u$ for every
  $u\geq0$ }. Remark that by monotonicity we can couple $s_t$ and $\hat
s_t$ in such a way that $s_t\subset\hat s_t$ for every $t\geq0$; also,
remark that once the deterministic set $S^-_t$ is empty there is no more constraint on the
dynamics (which does not mean that $s_t=\hat s_t$ after that time!).

\medskip

The basic point is the following:
\begin{Theorem}
\label{th:souffle'}
  W.l.p., for every $0\leq t\leq M^2$ one has
  \begin{eqnarray}
\label{SetaS}
S_t^-\subset    \hat s_t\subset S_t^+.
  \end{eqnarray}
\end{Theorem}

Of course, the lower bound is trivial by the very definition of $\hat
s_t$ (and holds not only w.l.p.\ but with probability one). 
As a side remark, it will be apparent below that it would be sufficent  to
have 
the above statement with $M^2$ replaced by $M^{1+\epsilon}$ for some $\epsilon>0$.

\bigskip

\noindent
{\sl Proof of Theorem \ref{th:ineqt} (assuming Theorem
  \ref{th:souffle'}).}
Thanks to \eqref{SetaS},
for all times 
$$(6/7) C(\alpha)M(\log M)^6\leq t\leq M^2$$ one has $\hat
s_t\subset S^+_t= B^-_{(0,0,0)}$, which is a subset of a cube of side
$(\log M)^2$. This is just because $S_t^+=B^-_{(0,0,0)}$ for $t\geq
(6/7) C(\alpha)M(\log M)^6$, as we observed just after the definition
of $S_t^+$ .
By point (1) in Lemma \ref{th:fnraf} below, this implies that, within
time $(6/7)C(\alpha)M(\log M)^6+O((\log M)^6)\leq C(\alpha)M(\log M)^6
\ll M^2$, $\hat s_t$ has hit the ground state $s^-$ at least once, with
probability at least $1-1/(2e)$. Since $s_t\subset \hat s_t$, this
implies that $s_t$ has also hit $s^-$ and \eqref{eq:ineqt} follows. \qed

\subsection{Controlling the auxiliary dynamics $\hat s_t$}

The first ingredient of the proof of Theorem \ref{th:souffle'} is
the following lemma, which gives a rough upper bound on the mixing time in a cube of size
$M$:
 \begin{Le}
 \label{th:fnraf}
 For every $\alpha>0$ there exists $C_1(\alpha)<\infty$
 such that for every $M\in \N$
   \begin{eqnarray}
     \label{eq:fnraf}
     \tmix\leq C_1(\alpha)M^3.
   \end{eqnarray}
   Moreover, there exists $C_2(\alpha)>0$ such that for every $T>0$
   the following is true:
\begin{enumerate}
\item the probability that $t(M)\geq T M^3$ is smaller than
$   \exp({-C_2(\alpha)T})$.
\item  with probability at least
$1-T\, e^{-C_2(\alpha)M }$,
the volume  of $s_t$ is at most
$M/10$ for all times $t(M)\leq t\leq T$.
\end{enumerate}
 \end{Le}
 \proof 
 We know that 
 \begin{eqnarray}
   \min_{\eta\in E_{\wedge,\vee}}\mu^{\wedge,\vee}(\eta)=
 \frac{e^{-2\alpha M^3}}{\hat Z_{\wedge,\vee}}\geq \frac{e^{-2\alpha M^3}}
 {\sum_{v\geq0}e^{-2\alpha v}\mathcal N(v)}\geq C_2(\alpha)e^{-2\alpha M^3},
 \end{eqnarray}
where $\mathcal N(v)$ was defined after formula \eqref{eq:Nv} and is just the number of
plane partitions of volume $v$.
 Then, it follows from \eqref{eq:tmixgap} (and the fact 
that the gap is uniformly positive) that the mixing time 
 is $O(M^3)$.

 From \eqref{eq:fnraf} it is immediate to deduce (modulo
 redefining $C_1(\alpha)$) that the probability that $t(M)\leq C_1(\alpha)
 M^3$ is greater than some $\epsilon(\alpha)>0$.  (Indeed, Lemma
 \ref{lemma_basico} implies that there exists a set $A$ of
 configurations such that $\mu^{\wedge,\vee}(A)>1/2$ and such that all
 the configurations $s\in A$ can be reached via at most $m=m(\alpha)$
 Markov Chain moves from the ground state $s^-$, for some
 $m(\alpha)<\infty$ independent of $M$).  From this, one easily
 deduces that, for every $n\in\N$, the probability that $t(M)> n\,
 C_1(\alpha)M^3$ is smaller than $(1-\epsilon(\alpha))^n$, i.e., claim
 (1).  Indeed, if the evolution has not hit $s^-$ before
 time $(n-1)\,C_1(\alpha)M^3$, just restart the dynamics from the
 maximal configuration $s^+$ at $t=(n-1)C_1(\alpha)M^3$: this can only
 make $t(M)$ larger, by monotonicity.

To prove statement (2), observe first of all that for all times $t>t(M)$ (or, more precisely,
conditionally on $t(M)<t$) the distribution of
$s_t$ is stochastically dominated by the equilibrium distribution $\mu^{\wedge,\vee}$.
On the other hand, Lemma \ref{lemma_basico} tells us that 
\begin{eqnarray}
\label{eq:equilib}
  \mu^{\wedge,\vee}(|s|\geq M/10)\leq \exp(-c_1(\alpha)M/10).
\end{eqnarray}
Secondly, the number of Markov Chain moves in the interval $[0,T]$ is a
Poisson random variable $\zeta$ with average $T M^2$, and an
 elementary computation shows that for a Poisson random variable $\zeta_\lambda$ of parameter $\lambda$ one has
\begin{eqnarray}
\label{eq:poiss}
  \bbP(\zeta_\lambda\geq n)\leq e^{-n(\log (n/\lambda)-1)}.
\end{eqnarray}
We have therefore, calling $t_i, i=1,\ldots,\zeta$ the random times when the updates occur,
\begin{eqnarray}
  \bbP\left(\exists t\in[t(M),T]:|s_t|>M/10\right)&\leq& e^{-M^2T}\\\nonumber
&&+
\bbP\left( \zeta\leq 4 M^2 T; \exists i: t(M)\leq t_i\leq T, |s_{t_i}|\geq M/10
\right)\\\nonumber
&&\leq  e^{-M^2T}+4 M^2\,T \,\exp(-c_1(\alpha)M/10),
\end{eqnarray}
where in the last inequality we used the union bound and  \eqref{eq:equilib}.
\qed

\medskip

{\sl Notational convention}: in the rest of this section, we will
use for simplicity of exposition expressions like ``for all times
larger than $t_0$'' to mean ``for all times $t_0\leq t\leq M^2$''.

\medskip

\begin{figure}[h]
\centerline{
\psfrag{x}{$x$}
\psfrag{y}{$y$}
\psfrag{z}{$z$}
\psfrag{s}{$\star$}
\psfrag{0}{$0$}
\psfig{file=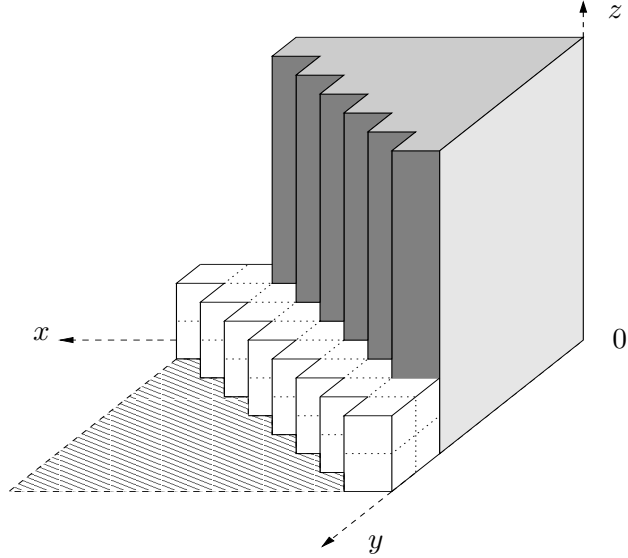,height=2.9in, width=3.2in}}
\caption{The maximal configuration $s^{\max}_{j}=:s^{(j)}_{t=0}$
  in $\mathcal S(j)$. The dark region is the one which is
  constrained to remain completely full during the evolution of
  $s^{(j)}_t$, while the dashed one remains empty.  The white
  region is the one which can evolve. Note that, for
graphical convenience, the axes are drawn with orientations which 
differ from those of Fig. \ref{fig:S-}.} 
\label{fig:Sjn}
\end{figure}

For the next lemma we need some notations.  Let $0\leq j<2(K-1)$ and
let $\mathcal S(j)$ be the set of configurations $s$ such that
$B_v\subset s$ if $(v_1+v_2)<j$ and $B_v\cap s=\emptyset$ if
$(v_1+v_2)>j+1$ or if $v_1+v_2\in\{j,j+1\}$ and $v_3\geq 2$. Let $
s^{\max}_{j}$ be the maximal configuration in $\mathcal S(j)$, 
see Figure \ref{fig:Sjn} (of course, both $\mathcal S(j)$ and $
s^{\max}_{j}$ depend on $M$).
We define a dynamics $\{s^{(j)}_t\}_{t\geq0}$ by requiring that
(in law) it equals our usual dynamics, with initial condition
$s^{(j)}_{t=0}=s^{\max}_{j}$ and conditioned on the event that
$s^{(j)}_t\in\mathcal S(j)$ for every $t\geq0$.

\begin{figure}[!h]
\centerline{
\psfrag{x}{$x$}
\psfrag{y}{$y$}
\psfrag{z}{$z$}
\psfrag{0}{$0$}
\psfrag{s}{$\star$}
\psfig{file=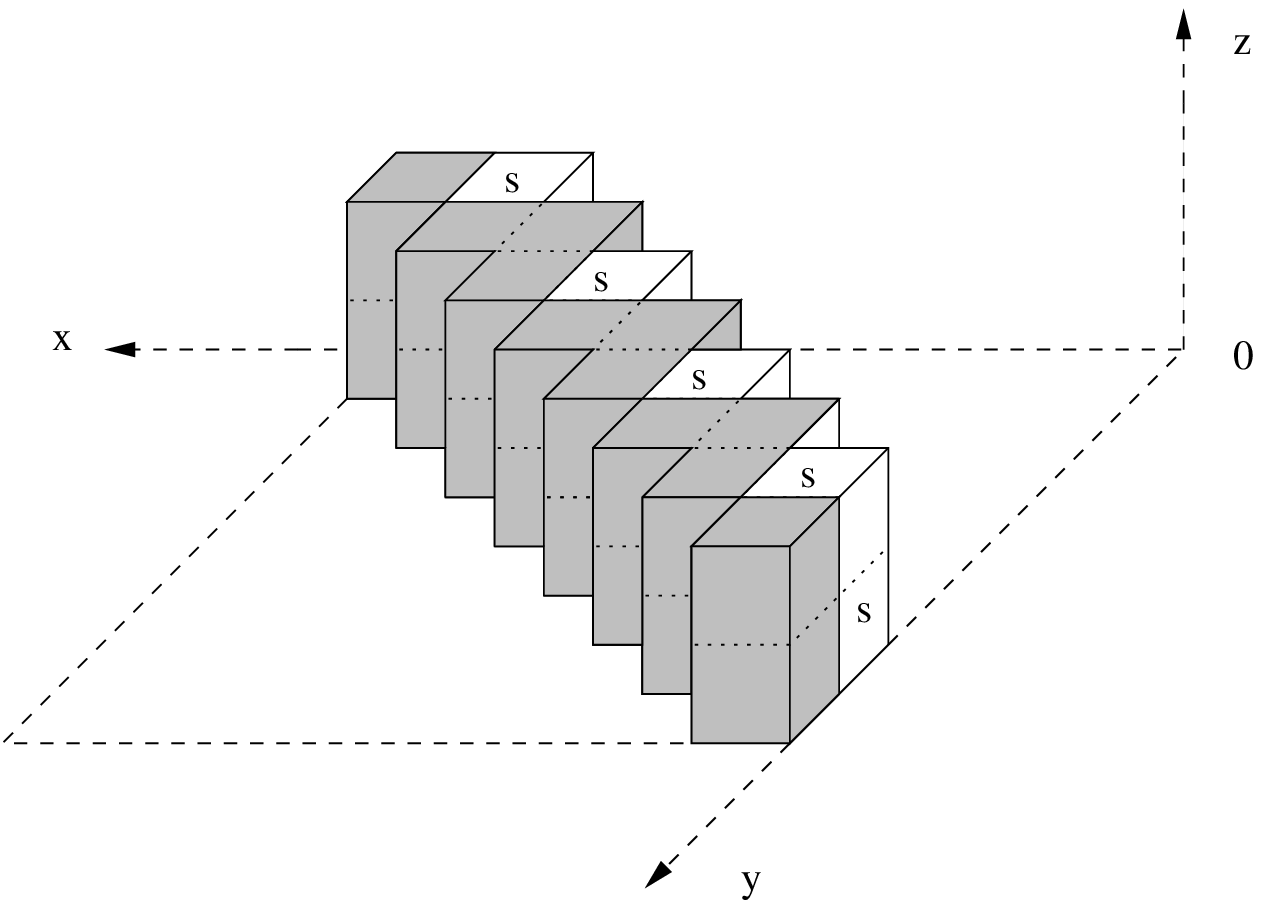,height=2.9in, width=3.2in}}
\caption{Here we drew only the region where $s^{(j)}_t$ is allowed
  to evolve (the white region in Fig.  \ref{fig:Sjn}). Under the
  dynamics $\hat s^{(j)}_t$, blocks marked by a $\star$ are
  constrained to remain full up to time $\tau/2$.  As a consequence,
  up to that time the remaining blocks (shown in gray) evolve
  independently by groups of at most six.}
\label{fig:Sjn2}
\end{figure}

\begin{Le}
\label{th:diag}
  W.l.p., the following holds for all times $t\geq \tau$. 
If $(v_1+v_2)=j$ and $v_3=1$, or if $(v_1+v_2)=j+1$, then 
$B_v\cap s^{(j)}_t=\emptyset$. If $(v_1+v_2)=j$ and $v_3=0$,
$B^+_v\cap s^{(j)}_t=\emptyset$.
\end{Le}
\proof
As in the proof of Theorem
\ref{th:ineqt}, we introduce an auxiliary dynamics $\{\hat
s^{(j)}_t\}_{t\geq0}$ for which the claim is easier to prove, and such
that almost surely $s^{(j)}_t\subset \hat s^{(j)}_t$ (so that the
claim follows also for $s^{(j)}_t$).  The auxiliary dynamics is
defined simply by requiring that its law equals that of
$\{s^{(j)}_t\}_{t\geq0}$ conditioned on the event that, for all times
$0\leq t\leq\tau/2$, $B_v\subset s^{(j)}_t$ if $(v_1+v_2)=j$ and $v_1\in
2\N+1$. In other words, such blocks $B_v$ are frozen and remain
completely full up to time $\tau/2$.  This implies that during this
time interval the blocks which are not frozen evolve independently by
groups of at most six, see Figure \ref{fig:Sjn2}.
 Applying Lemma \ref{th:fnraf}, one sees that at time $\tau/2$ the situation 
is the following (w.l.p.):
\begin{enumerate}
\item $B_v\cap\hat s^{(j)}_{\tau/2}=\emptyset$ if $(v_1+v_2)=j+1$ or
  $(v_1+v_2)=j,v_1\in 2\N$ and $v_3=1$; the same remains true for
  all later times (of course here we are using the fact that, since $\tau=
  C(\alpha)(\log M)^8/7$, 
 the probability that $t((\log M)^2)\geq\tau/2$ is
$O(\exp(-C'(\alpha)(\log M)^2))$, cf. statement (1) of Lemma \ref{th:fnraf}, 
i.e., w.l.p.\ such event does not occur.)
\item If $(v_1+v_2)=j,v_1\in 2\N$ and $v_3=0$, the volume of
 $B_v\cap\hat s^{(j)}_{\tau/2}$ is at most $(\log M)^2/10$ (in particular,
$B_v^+\cap \hat s^{(j)}_{\tau/2}=\emptyset$); the
  same remains true for all later times. Here we are using statement (2) of Lemma \ref{th:fnraf},
with $T=M^2$.
\item all blocks marked by a star in Figure \ref{fig:Sjn2} are of course 
still completely full.
\end{enumerate}

During the interval $(\tau/2,\tau]$ the blocks marked by a $\star$ in
Figure \ref{fig:Sjn2} are free to move.  However, now the pairs of
blocks $v_1+v_2=j,v_1\in 2\N+1, v_3\in \{0,1\}$ evolve independently for
$v_1$ different, thanks to point (1) above.  Therefore, applying once
more Lemma \ref{th:fnraf} and point (2) above, one sees that for
$t\geq\tau$ one has $B_v\cap\hat s^{(j)}_{t}=\emptyset$ for
$v_1+v_2=j,v_1\in 2\N+1, v_3=1$, and $B_v^+\cap\hat
s^{(j)}_{t}=\emptyset$ for $v_1+v_2=j,v_1\in 2\N+1, v_3=0$. The claim is
proven.  \qed

\bigskip

\noindent
{\sl Proof of Theorem \ref{th:souffle'}.} It is sufficient to
prove the following: for every $i\in\N\cup\{0\}$ such that $i\tau\leq
M^2$, one has (w.l.p.)
\begin{eqnarray}
  \label{eq:induz2}
  \hat s_t\subset S^+_{i\tau }\;\;\mbox{for every}\;\; t\geq i\tau.
\end{eqnarray}

For $i=0$ the statement 
is trivial since $S_0^+=\mathcal C_M$.

Now we assume that the claim is true up to a certain $i$, and we
show that it holds also for $i+1$.  It is convenient to introduce the
definition of column $C_{v_1,v_2}(H)$ with (integer) base coordinates $0\leq
v_1,v_2<K$ and height $0\leq H\leq K$: this is just a parallelepiped of
height $(\log M)^2 H$ whose base is the square of side $(\log M)^2$
 such that the base point of minimal $L^1$ norm has coordinates
$((\log M)^2v_1,(\log M)^2 v_2)$. Note that, for each $t$, $S^\pm_t$ can be
viewed as composed of $K^2$ such columns: in the case of $S_t^-$ the
heights $H_t^-(v_1,v_2)$ take values in $\{0,1,\ldots,K\}$ while in
the case of $S_t^+$ we call them $H_t^+(v_1,v_2)$ and they take
values in $\{0,1/2,1,\ldots,K-1/2,K\}$.  Note also that, by
construction, $H_t^\pm(v_1,v_2)$ depends on $(v_1,v_2)$ only through
$v_1+v_2$ (see Figure \ref{fig:S-}).

To complete the inductive proof, we need to prove that for all $(v_1,v_2)$
\begin{eqnarray}
  \label{eq:induc2}
 \left[ \hat s_t\cap  C_{v_1,v_2}(K)\right]\subset C_{v_1,v_2}(H^+_{(i+1)\tau}
(v_1,v_2))
\;\;\mbox{for every}\;\; t\geq (i+1)\tau.
\end{eqnarray}
The following cases can occur (keep Figure \ref{fig:S-} in mind):
\begin{enumerate}
\item 
$H_{(i+1)\tau}^-(v_1,v_2)=K$. In this case, \eqref{eq:induc2} is obvious
because also $H_{(i+1)\tau}^+(v_1,v_2)=K$.
\item 
$1<H_{(i+1)\tau}^-(v_1,v_2)<K$. In this case, 
we have that 
\begin{eqnarray}
\label{eq:414}
 C_{v_1,v_2}(H^-_{(i+1)\tau}(v_1,v_2))\subset
   \left[ \hat s_t\cap  C_{v_1,v_2}(K)\right]
\subset C_{v_1,v_2}(H^-_{(i+1)\tau}(v_1,v_2)+2)
\end{eqnarray}
for every $i\tau<t<(i+1)\tau$. The lower bound is trivial by the
definition of the dynamics $\hat s_t$, while the upper bound follows
from the inductive hypothesis \eqref{eq:induz2} and from the fact that
the definition of $S_t^+$ implies that $0\leq
H^+_t(v_1,v_2)-H_t^-(v_1,v_2)< 2$.  
Since we want to prove  \eqref{eq:induc2}, by monotonicity we can assume that all the columns
labeled $(w_1,w_2)$ with $w_1+w_2<v_1+v_2$ are completely full during the time  interval
$i\tau<t<(i+1)\tau$.
But then, as we shall argue in a moment, in the time
interval $(i\tau,(i+1)\tau)$ the column $\hat s_t\cap C_{v_1,v_2}(K)$
evolves independently of all the others, and an application of
Lemma \ref{th:fnraf} implies \eqref{eq:induc2}, since
$H^+_{(i+1)\tau}(v_1,v_2)$ is just $H^-_{(i+1)\tau}(v_1,v_2)+1/2$,
cf.\ the caption of Figure \ref{fig:S-}. 

To see that the column $(v_1,v_2)$ evolves 
independently of all the others in the interval $i\tau<t<(i+1)\tau$, note that 
it can be influenced only by the columns labeled  $(w_1,w_2)$ with $w_1+w_2=v_1+v_2+1$.
However (cf.\ Figure \ref{fig:S-}) in this case $H^-_{(i+1)\tau}(w_1,w_2)=H_{(i+1)\tau}^-(v_1,v_2)-2$
so that, from the induction hypothesis (cf. \eqref{eq:414})
\begin{eqnarray}
\hat s_t\cap C_{w_1,w_2}(K)\subset C_{w_1,w_2}(H_{(i+1)\tau}^-(v_1,v_2)).
\end{eqnarray}
In other words, the column $(w_1,w_2)$ is too low to influence the column $(v_1,v_2)$.

\item $H_{(i+1)\tau}^-(v_1,v_2)=1$. Again one has \eqref{eq:414} and one can assume by
monotonicity that the columns with $w_1+w_2<v_1+v_2$ are completely full in the time interval
$i\tau<t<(i+1)\tau$. The argument proceeds like in the previous case once one realizes that 
for $w_1+w_2=v_1+v_2+1$ one has $H^+_{i\tau}(w_1,w_2)=1/2$, so that by the induction 
hypothesis $\hat s_t\cap C_{w_1,w_2}(K)\subset C_{w_1,w_2}(1/2)$ for $t\ge i\tau$ and such column 
cannot influence the one labeled $(v_1,v_2)$.

\item It remains to consider the case of the columns with $H_{(i+1)\tau}^-=0.$ 
Define
$$
j:=\max_{0\le w_1<K,0\le w_2<K} \{ w_1+w_2: \; H_{(i+1)\tau}^-(w_1,w_2)>0\}+1 < 2K-1,
$$
with the convention that $j:=0$ if the set is empty.
It is convenient to distinguish two sub-cases:
\begin{enumerate}
\item if $v_1+v_2>j+1$ then, by definition of
  $S_t^+$ one sees that $H_{i\tau}^+(v_1,v_2)=0$, so that
  \eqref{eq:induc2} follows (both sets are empty) from the induction
  hypothesis \eqref{eq:induz2}. 
\item  if  $j\le v_1+v_2\le j+1$ then by monotonicity we can assume that all columns with 
$w_1+w_2< j$ are completely full in the time interval $i\tau\le t\le (i+1)\tau$, and on the other hand 
we know that all columns with $w_1+w_2>j+1$ are
completely empty for $t\ge i\tau$. Also, we know
from \eqref{eq:induz2} that, always for $t\geq i\tau$,  $[\hat s_t\cap
C_{w_1,w_2}(K)]\subset C_{w_1,w_2}(2)$ for all $(w_1,w_2)$ such that $w_1+w_2\in\{j,j+1\}$.
We can therefore apply Lemma \ref{th:diag} to deduce that for all times
$t\geq(i+1)\tau$ the columns with $w_1+w_2=j+1$ are completely empty,
while if  $w_1+w_2=j$ then $[\hat s_t\cap
C_{w_1,w_2}(K)]\subset B^-{(w_1,w_2,0)}$. Recalling the definition of
$S^+_t$ (in particular, as explained in the caption of Figure \ref{fig:S-})
we have therefore
proven \eqref{eq:induc2} for all columns $(w_1,w_2)$ such that  $w_1+w_2\in
\{j,j+1\}$. 
\end{enumerate}
\end{enumerate}

\qed

\subsection{The general case}

\label{sec:xisigma}
Here we prove Theorem \ref{th:tmix1} in the general situation where
$\xi$ is not the maximal configuration in $\Omega^h_{1,n}$.  The proof
is conceptually similar to the one where the ceiling is maximal, and
therefore some arguments will be only sketched (see however
 Remark \ref{rem:frastagl} below, where an important difference
between the two cases is pointed out).

As in section \ref{sec:maxceil}, we can assume by monotonicity that
$h=0$, $n=2M$, $k=M$ and $\si=\vee$.  First of all, it is important to
realize that the maximal configuration (ground state)
$\{\eta^{(j)}=\xi$ for all $j=1,\ldots,M\}$ corresponds to the subset
of $\mathcal C_M$ defined by the property that the elementary cube
labeled $r=(r_1,r_2,r_3)$, with $r_i\in \{0,\ldots,M-1\}$ (cf. Section \ref{sec:polyset})
belongs to $s^-$ if and only if
\begin{eqnarray}
  \xi_{M-r_1+r_2}<M-r_1-r_2
\end{eqnarray}
(we still call the ground state $s^-$, even if it is no longer the empty
set as in previous section).  We note also that the equilibrium
measure \eqref{eq:moreexp} is given in this case by
\begin{eqnarray}
  \mu^{\xi,\vee}(s)=\frac{e^{-2\alpha|s\setminus s^-|}
  }{\sum_{s'\supset
s^-}e^{-2\alpha|s'\setminus s^-|}
}{\bf 1}_{s'\supset s^-}.
\end{eqnarray}
Of course, the sum in the denominator is only over the configurations
$s'\subset \mathcal C_M$ which are given by unions of elementary cubes.
The dynamics $s_t$ just coincides with that described in Section
\ref{sec:polyset}, {\sl conditionally on the event} that $s_t\supset s^-$ for
every $t$ (of course, the initial condition $s_{t=0}$ has to
verify the same property).

\begin{figure}[h]
\centerline{
\psfrag{x}{$x$}
\psfrag{y}{$y$}
\psfrag{0}{$0$}
\psfrag{s}{$\star$}
\psfig{file=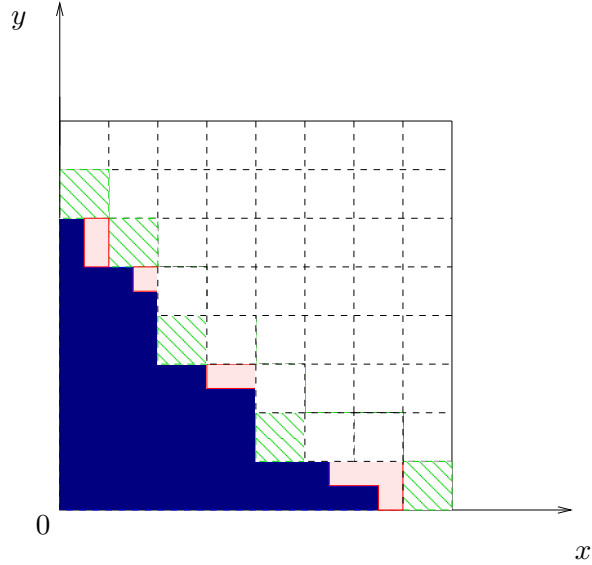, width=3in}}
\caption{The sets $A_{i}$ and $s^-$ seen from above. Squares should be
  imagined to have side $(\log M)^2$, and here $K=M/(\log M)^2=8$.
  The dark region is a horizontal section of $s^-$, while the lightly
  colored (respectively, the dashed) region is the set $A_2\setminus
  s^-$ (resp.\ $A_3$) seen from above. $s^-$ extends vertically up to
  height $M$, while $A_2,A_3$ extend  only up to height
  $(\log M)^2$ (i.e., one single block).  Note that the block
  $B_{(0,6,0)}$ belongs to $A_3$ because neither $B_{(-1,6,0)}$ nor
  $B_{(0,5,0)}$ belong to $\mathcal C_M\setminus A_2$ (in particular,
  $B_{(-1,6,0)}$ falls out of $\mathcal C_M$). A similar remark holds
  for  $B_{(7,0,0)}$.  }
\label{fig:A}
\end{figure}

\medskip
As in Section \ref{sec:proofof}, we cut $\mathcal C_M$ into $K^3$ blocks
$B_v$ of side $(\log M)^2\in\N$. 
Consider the dynamics started from some $s_{t=0}\supset s^-$. The basic
estimate which allows to prove Theorem \ref{th:tmix1} is the following:
\begin{Proposition}
\label{th:propA3}
For every $\alpha>0$ there exists 
$C(\alpha)<\infty$ (independent of $M$, $\xi$ and $s_{t=0}$) such that
 the following holds w.l.p.: 
 \begin{eqnarray}
   s_t\subset (s^-\cup A_1)
 \end{eqnarray}
for every $(6/7)C(\alpha)M(\log M)^6\leq t\leq M^2$
where $A_1=A_2\cup A_3$,
\begin{eqnarray}
  A_2:=\cup \{B_v: v_3=0, B_v\cap s^- \mbox{\;contains at least one
elementary cube} \}
\end{eqnarray}
and 
\begin{eqnarray}
  A_3:=\cup\{B_v:
  v_3=0, B_{(v_1-1,v_2,0)}\notin (\mathcal C_M\setminus A_2) \mbox{\;and\;}B_{(v_1,v_2-1,0)}\notin 
(\mathcal C_M\setminus A_2)\},
\end{eqnarray}
see Figure \ref{fig:A}.
\end{Proposition}
Note that, in the case $\xi=\wedge$, one has $A_2=\emptyset$, $A_3=B_{(0,0,0)}$  and Proposition 
 \ref{th:propA3} follows from Theorem \ref{th:souffle'} above.

Proposition \ref{th:propA3} is proven below, and now we show that it does
imply Theorem \ref{th:tmix1}. 

\smallskip

{\sl Proof of Theorem \ref{th:tmix1} (for general ceiling $\xi$) assuming
Proposition \ref{th:propA3}.}

Note that the volume of $A_1$ satisfies
\begin{eqnarray}
  |A_1\setminus s^-|\leq 4 M (\log M)^4,
\end{eqnarray}
so that 
\begin{eqnarray}
\label{eq:mins}
\min_{s:s^-\subset s\subset
 (s^-\cup A_1)}
 \mu^{\xi,\vee}\left(\left.s\right|  s^-\subset s\subset
 (s^-\cup A_1)\right)\geq e^{-c_8(\alpha)M(\log M)^4}
\end{eqnarray}
for some $c_8>0$ (we used the fact that the number of configurations $s$ satisfying 
$s^-\subset s\subset
 (s^-\cup A_1)$ is smaller than 
$$
2^{|A_1\setminus s^-|}:
$$
this would be the exact number of configuration if there were no monotonicity constraints on $s$).
Call, for ease of notation,
\begin{eqnarray}
  T_M:=C(\alpha) M(\log M)^6.
\end{eqnarray}
From  \eqref{eq:tmixgap}, \eqref{eq:mins} and Theorem
\ref{teo_gap} one easily deduces that the mixing time of the dynamics
constrained to $s^-\subset s_t\subset [s^-\cup A_1]$
is $O(M(\log M)^4)$. Since $(6/7)T_M+M(\log M)^4\ll
T_M\ll M^2$, Theorem \ref{th:tmix1} easily follows. 
\qed

\begin{remark}\rm
\label{rem:frastagl}
 It is important to notice a crucial difference between the proof
 of   Theorem \ref{th:tmix1} for maximal ceiling $\wedge$ and for
 generic ceiling $\xi$. In the former case, we have shown a stronger 
statement, i.e. (cf. Theorem 
\ref{th:ineqt}) that  the
dynamics started from an arbitrary $s_{t=0}$ hits the ground state
$s^-$ within a time of order $M(\log M)^6$. This is related to the fact that the ground state $s^-$ has an
equilibrium weight which is positive, uniformly in the system size. 
In the generic case, i.e. when 
 the ceiling $\xi$ has a jagged shape, there is no reason to believe
 that
the dynamics hits $s^-$ within such time;  indeed, its equilibrium
 weight (an therefore the inverse of its hitting time) may well be exponentially small in  $M$.
\end{remark}

\subsubsection{Proof of Proposition \ref{th:propA3}}
By monotonicity, it is clear that it suffices to prove the claim for 
$s_{t=0}=\mathcal C_M$.
As in the proof of Theorem \ref{th:ineqt}, we introduce two deterministic
subsets $\hat S_t^\pm$ of $\mathcal C_M$. If $S_t^\pm$ are the sets
which were defined in Section \ref{sec:proofof}, then we establish that
\begin{eqnarray}
  \hat S_t^-:=S_t^-\cup s^-,
\end{eqnarray}
while 
\begin{eqnarray}\hat S_t^+:= S_t^+\cup A_1.
\end{eqnarray}
We note that from the discussion of the 
properties of  $S_t^\pm$ in Section \ref{sec:proofof} it follows that
 for $t\geq (6/7) T_M$ one has $\hat S_t^-=s^-$ and $\hat S_t^+= A_1$. The claim of the proposition then follows if we can 
prove, in analogy with Theorem \ref{th:souffle'}, that
w.l.p.\ one has 
\begin{eqnarray}
\label{eq:424}
  s_t\subset \hat S_t^+
\end{eqnarray}
for every $0\leq t\leq M^2$.  The proof of this fact proceeds with the help of an
auxiliary dynamics $\hat s_t$, which dominates $s_t$, and whose law is
that of $s_t$ conditioned on the event that $\{s_u\supset \hat S_u^-$
for every $u\geq0\}$. Since  the proof is very similar to
that of Theorem \ref{th:souffle'}, we do not give details.
The only fact which requires some care is that we cannot
apply Lemma \ref{th:fnraf} to the blocks which have a non-empty intersection
with $s^-$, since a certain number of its elementary cubes are frozen to be full for all times.

\smallskip

The extra result we need concerns therefore the evolution of a {\sl
  single cube} of side $(\log M)^2$ and with arbitrary ceiling. Let $\tilde \xi\in \Omega_{1,2(\log
  M)^2}^0$, let $\tilde s^-\subset\mathcal C_{(\log M)^2}=B_{(0,0,0)}$ be the
ground state corresponding to the ceiling $\tilde \xi$ and $\tilde s_t$ be
the evolution started from the full configuration $B_{(0,0,0)}$ and
constrained to
$$
\tilde s^-\subset \tilde s_t\subset B_{(0,0,0)}
$$ for
every $t>0$, and call $\tilde T_{mix}$ its mixing time. Needless to say,
 its invariant measure is $\mu^{\tilde\xi,\vee}$. Notice that
when $\tilde s^-=B_{(0,0,0)}$ (i.e., when $\tilde \xi=\vee$) the dynamics is trivial
($\tilde s_t=B_{(0,0,0)}$ for all times) while when $\tilde s^-=\emptyset$ (i.e., when
$\tilde\xi=\wedge$) 
the forthcoming lemma is already implied by Lemma \ref{th:fnraf} (just replace $M$ with $(\log M)^2$ there).
 \begin{Le}
 \label{th:fnraf2}
 For every $\alpha>0$ there exist $C_4(\alpha)<\infty$ and $C_5(\alpha)>0$
 such that uniformly in  $M$ and $ \tilde\xi$ 
   \begin{eqnarray}
     \label{eq:fnraf2}
     \tilde T_{mix}\leq C_4(\alpha)(\log M)^6.
   \end{eqnarray}
Moreover, with $\tau$  defined as in \eqref{eq:tau},
\begin{eqnarray}
\label{eq:munu}
  \bbP\left(\exists\; \tau\leq t\leq M^2: (\tilde s_t\setminus \tilde s^-)\cap B_{(0,0,0)}^+\ne \emptyset\right)\leq
e^{-C_5(\alpha)(\log M)^2}.
\end{eqnarray}
 \end{Le}
What Eq.
 \eqref{eq:munu} is saying is essentially that for all times larger
 than $(\log M)^2\tilde T_{mix}$ but smaller than $M^2$ the upper half
 of the cube under consideration contains only the elementary cubes
 which are imposed by the constraint $\tilde s^-\subset
 \tilde s_t$.  Useless to say, this means that w.l.p.\ the event in the
 left-hand side of \eqref{eq:munu} does not occur.

Using Lemma \ref{th:fnraf2} as a substitute for Lemma \ref{th:fnraf},
the proof of Proposition \ref{th:propA3} is easily concluded in
analogy with the proof of Theorem \ref{th:souffle'}. We spare the
reader additional details.
\qed

\medskip

{\sl Proof of Lemma \ref{th:fnraf2}.}
The proof of \eqref{eq:fnraf2} is completely analogous to that of \eqref{eq:fnraf}
and uses the fact that the spectral gap of the dynamics is positive, uniformly in 
$M$ and $\tilde\xi$.
To prove \eqref{eq:munu}, let us recall the well known inequality
which relates the total variation distance from equilibrium
of a reversible Markov Chain at time $t$ with its mixing time $\tmix$:
\begin{eqnarray}
\label{eq:wk}
  \sup_{s\in\Omega}|| P_t^s(\cdot)-\mu(\cdot)||_{\rm var}\leq e^{-\lfloor t/\tmix\rfloor},
\end{eqnarray}
where $\Omega$ is the state space of the Markov Chain, $\mu$ its
invariant measure and $P_t^s$ the law at time $t$, if the initial
condition at time zero is $s$.  Call $t_i, 1\leq i\leq \zeta$ the random
times when Markov Chain updates occur in the time interval
$[\tau,M^2]$, and observe that in our case $\zeta$ is a Poisson random
variable of parameter $(M^2-\tau)(\log M)^4$.  One has then, using
\eqref{eq:poiss}, that the left-hand side of \eqref{eq:munu} is upper
bounded by
\begin{eqnarray}
\label{eq:428}
  e^{-M^2}+\bbP\left[\zeta\leq 4 M^2(\log M)^4\mbox{\;and\;}\exists \; i\leq \zeta: 
(\tilde s_{t_i}\setminus \tilde s^-)\cap B_{(0,0,0)}^+\ne \emptyset\right].
\end{eqnarray}
On the other hand, calling $\tilde P_t(\cdot)$ the law of $\tilde s_t$ and defining the set
$$U:=\{s\subset B_{(0,0,0)}: (s\setminus \tilde s^-)\cap B_{(0,0,0)}^+\ne \emptyset\},$$ one has
\begin{eqnarray}
  \bbP\left[(\tilde s_t\setminus s^-)\cap B_{(0,0,0)}^+\ne \emptyset\right]&=&
  \tilde P_t\left(U\right) \leq\left|\tilde P_t\left(U\right)-\mu^{\tilde\xi,\vee}\left(U\right) 
\right|+\mu^{\tilde\xi,\vee}\left(U\right)\\\nonumber
&\leq & e^{-\lfloor t/\tilde T_{mix}\rfloor}+e^{-C_5(\alpha)(\log M)^2}.
\end{eqnarray}
In the last step, we used \eqref{eq:wk} for the first term and Lemma \ref{lemma_spiaccicato}
to estimate the equilibrium probability of $U$. 
The claim \eqref{eq:munu} then follows from \eqref{eq:428}, a union bound and the
fact that $\tau/\tilde T_{mix}\geq C_6(\alpha)(\log M)^2$ for some $C_6(\alpha)>0$.  \qed

\bigskip
%

\section*{Acknowledgments} F.L.T.\ was partially supported by 
ANR, project POLINTBIO and project LHMSHE. P.C.\ was partially supported by 
NSF Grant DMS-0301795. F.M. was partially supported by the Advanced Research Grant ``PTRELSS''
ADG-228032 of the European Research Council.


\begin{thebibliography}{99}
\bibitem{BD} R. Bubley and M. Dyer, {\sl Path coupling: A technique
    for proving rapid mixing in Markov chains}, Proc. of the 38th
    Annual Symposium on Foundations of Computer Science (1997), 223--231. 


\bibitem{CK} R.\ Cerf, R.\ Kenyon, {\em The low-temperature expansion
    of the Wulff crystal in the 3D Ising model}, Comm.\ Math.\ Phys.\
  {\bf 222} (2001), 147--179. MR1853867


\bibitem{CLP} H.\ Cohn, M.\ Larsen, J.\ Propp, {\em The shape of a
    typical boxed plane partition}, New York J.\ Math.\ {\bf 4} (1998),
  137--165. MR1641839 

\bibitem{CM} P.\ Caputo, F.\ Martinelli, {\em 
Asymmetric diffusion and the energy gap above the 111 ground state of the quantum XXZ model},
Comm.\ Math.\ Phys.\ {\bf 226} (2002), 323--375. MR1892457

\bibitem{CMT} P.\ Caputo, F.\ Martinelli, F.L.\ Toninelli, {\em On the
    approach to equilibrium for a polymer with adsorption and
    repulsion}, Electronic Journal of Probability {\bf 13} (2008), 213--258. MR2386733

\bibitem{GPR} S.\ Greenberg, A.\ Pascoe, D.\ Randall, {\em Sampling
    biased lattice configurations using exponential metrics}, 
  Proc. of the 19th Annual ACM-SIAM Symposium on
  Discrete Algorithms (2009),76--85.

\bibitem{K} R.\ Kenyon, {\em Lectures on dimers}, 
Statistical mechanics, 
IAS/Park City Math.\ Ser.\ {\bf 16} (2009), 191--230. 
Amer. Math. Soc., Providence.  MR2523460

\bibitem{cf:peres-book} D. A. Levin, Y. Peres, E. L. Wilmer, 
{\sl Markov Chains and Mixing Times}, American Mathematical Society (2009). MR2466937

\bibitem{Wilson}  D.\ B.\ Wilson, {\em Mixing times of Lozenge
    tiling and card shuffling Markov chains},  Ann.\ Appl.\ Probab.\
  {\bf 14} (2004), 274--325.  MR2023023


\bibitem{cf:Wright} E. M. Wright, {\sl Asymptotic partition formula, 
I. Plane partitions}, Quart. J. Math. {\bf 2} (1931), 177-189.
\end{thebibliography}
\end{document}